\documentclass[11pt,a4paper]{amsart}

\usepackage[]{graphicx}
\usepackage{psfrag}

\newcommand{\ncmd}{\newcommand}
\ncmd{\rencmd}{\renewcommand}
\ncmd{\dspst}{\displaystyle }

\ncmd{\preu}{\noindent \mbox{\bf Proof  }}
\ncmd{\sketch}{\noindent \mbox{\bf Sketch of proof  }}
\ncmd{\preud}[1]
{\noindent \mbox{\bf Proof of \ref{#1}   }}

\ncmd{\ex}{\noindent \mbox{\bf exemple  }}
\ncmd{\exs}{\noindent \mbox{\bf exemples  }}
\ncmd{\fin}{\hspace*{\fill} 
\quad\hbox{\hskip 1pt\vrule width 4pt height 6pt
          depth 1.5pt\hskip 1pt} \medskip }
\ncmd{\lemv}{\noindent \mbox{\bf Lemma  }}
\ncmd{\thmv}{\noindent \mbox{\bf Theorem  }}

\ncmd{\ab}{\mbox{$\mbox{AdS}_{n}\/$ }}
\ncmd{\ei}{\mbox{${\mbox{Ein}}_{n}\/$ }}
\ncmd{\hei}{\mbox{$\widehat{\mbox{Ein}}_{n}\/$ }}
\ncmd{\tei}{\mbox{$\widetilde{\mbox{Ein}}_{n}\/$ }}

\ncmd{\mr}{\mbox{$\widetilde{M}\/$ }}

\ncmd{\slrr}{\mbox{$\widetilde{SL}(2,{\mathbb R}) \times \widetilde{SL}(2,{\mathbb R})\/$ }}

\newtheorem{prop}{Proposition}[section]
\newtheorem{thm}[prop]{Theorem}
\newtheorem{lem}[prop]{Lemma}
\newtheorem{cor}[prop]{Corollary}
\newtheorem{defin}[prop]{Definition}
\newtheorem{rac}[prop]{Remark}
\newtheorem{ric}[prop]{Example}

\ncmd{\exe}{\begin{ric} \em }
\ncmd{\eexe}{\em \end{ric}}
\ncmd{\rque}{\begin{rac} \em}
\ncmd{\erque}{\em \end{rac}}

\rencmd{\thesubsubsection}{\thesubsection-\alph{subsubsection}}
\newcounter{inc}
\rencmd{\theequation}{\arabic{inc}}

\makeindex

\begin{document}

\vspace{5cm}

\title[BTZ multi black-holes]{Causal properties of AdS-isometry groups II: \\
BTZ multi black-holes}
\author{Thierry Barbot}
\email{Thierry.Barbot@umpa.ens-lyon.fr}
\thanks
{Work supported 
by CNRS, ACI ``Structures g\'eom\'etriques et Trous Noirs''.}
\address{CNRS, UMR 5669\\
Ecole Normale Sup\'erieure de Lyon\\ 
46 all\'ee d'Italie, 69364 Lyon}
\date{\today}

\begin{abstract}
This paper is the continuation of \cite{ba1}. We essentially prove
that the familly of strongly causal spacetimes defined in \cite{ba1}
associated to generic achronal subsets in $\mbox{Ein}_2$ contains 
all the examples of BTZ multi black-holes. It provides new elements
for the global description of these multi black-holes.
We also prove that any strongly causal spacetime locally modelled on 
the anti-de Sitter
space admits a well-defined maximal strongly causal
conformal boundary locally modelled
on $\mbox{Ein}_2$.
\end{abstract}

\maketitle

\section{Introduction}

In \cite{ba1}, we studied certain aspects of causal properties of 
AdS-spacetimes, i.e.
lorentzian manifolds of dimension $3$ with negative constant curvature; 
in other words,
locally modeled on the anti-de Sitter space AdS. In particular, we made
a detailed analysis on the causality notion in AdS, which, to be meaningfull,
has to be understood as the projection of the causality relation 
in the universal covering $\widetilde{\mbox{AdS}}$. The Einstein universes
$\mbox{Ein}_3$ and $\mbox{Ein}_2$ play an important role in this study: the
first one as containing a conformal copy of AdS, the second being the
conformal boundary of AdS (see {\S} $4$ in \cite{ba1}).

We also introduced the
notion of \emph{(generic) closed achronal subset\/} of the conformal 
boundary $\mbox{Ein}_2$ ({\S} $5$ of \cite{ba1}). We associate to it 
the invisible domains 
$E(\Lambda) \subset \mbox{AdS}$ and $\Omega(\Lambda) \subset \mbox{Ein}_2$ 
({\S} $8$). One of the main results of \cite{ba1} was the following:
assume that $\Lambda$ is non-elementary, and that a torsion-free
discrete group of isometries of AdS preserves $\Lambda$. Then, the
action of $\Gamma$ on $E(\Lambda)$ is free, properly discontinuous, and 
the quotient space $M_\Lambda(\Gamma) = \Gamma\backslash{E}(\Lambda)$ is
strongly causal (Theorem $10.1$). We observed that this construction
provides, when $\Lambda$ is a topological circle, the entire family 
of maximal globally hyperbolic AdS-spacetimes
admitting a Cauchy-complete Cauchy surface ({\S} $11$; the same
result is obtained in \cite{BenBon2} with a different proof). 
We also discussed
in which condition a discrete group of isometries $\Gamma$ is admissible, i.e.
preserves a closed achronal subset as above: essentially, up to a permutation
of space and time, admissible groups are precisely the positively proximal 
groups which are precisely the groups preserving a proper closed convex domain
of ${\mathbb R}P^3$ (Proposition $10.23$ of \cite{ba1}). We also stressed
out that there is a natural 1-1 correspondance between admissible groups
and pairs of (marked) complete hyperbolic metrics on the same surface
(not necesserily closed).
Finally, if $\Gamma$ is admissible and non-abelian, there exists 
an unique minimal closed generic achronal $\Gamma$-invariant subset 
$\Lambda(\Gamma)$ which is contained in every closed achronal 
$\Gamma$-invariant subset (Theorem~$10.13$, 
corollary~$10.14$ of \cite{ba1}).
Hence, following the classical terminology used for isometry 
groups of ${\mathbb H}^n$, it is natural to call $\Lambda(\Gamma)$ 
the \emph{limit set} of $\Gamma$.

The present paper completes this study in the elementary case, i.e. 
the case where $\Lambda$ is contained in the past or future of one point 
in $\mbox{Ein}_2$. In particular, we give the description of invisible
domains from elementary generic achronal domains (\S~\ref{sub.elementary}).

There is another equivalent definition of invisible domains 
$E(\Lambda(\Gamma))$ when  $\Lambda(\Gamma)$ is 
the limit set of $\Gamma$ (since $\Lambda(\Gamma)$ is minimal, 
$E(\Lambda(\Gamma))$ is in some way maximal among the $\Gamma$-invariant 
invisible domains). For any element $\gamma$ of $\Gamma$ we define the 
\emph{standard causal domain} $C(\gamma)$ as the set of points $x$ 
in AdS which are not causally 
related to their image $\gamma{x}$. 
The interior of the intersection of all the $C(\gamma)$ 
is then the standard causal domain of $\Gamma$, denoted 
$C(\Gamma)$. Every element $\gamma$ is also the time $1$ 
of a Killing vector field. We then define the 
\emph{absolute causal domain} $D(\gamma)$ as 
the domain of AdS where this Killing vector field 
has positive norm. The interior of the intersection of 
all the $D(\gamma)$ is the absolute causal domain $D(\Gamma)$ 
of $\Gamma$. 
According to Theorem~\ref{thm.egal} and corollary~\ref{cor.identifionsC}, 
\emph{if $\Gamma$ is a non-cyclic admissible group, then 
$D(\Gamma)$ and $C(\Gamma)$ both coincide with 
$E(\Lambda(\Gamma))$.\/} A first version of this theorem, 
in a very particular case, appears in \cite{amispin}, {\S}~$7$.

But, most of
all, we prove that the spacetimes $M_\Lambda(\Gamma)$ form a global family
of \emph{BTZ multi black-holes,\/} containing and enlarging all the
previous examples 
(\cite{amaking, aminneborg, amispin, BTZ, BTZ2, brill1, brill}),
except the non-static single BTZ black-holes themselves 
(see remark~\ref{rk.singleBTZ})

In order to develop this assertion, since we don't assume any acquaintance
of the reader with this notion, we need some preliminary discussion 
about black-holes.

\subsection{A quick insight into Schwarzschild and Kerr black-holes}
A honest introduction to the notion of black-hole requires a minimal 
historical exposition. First of all, a black-hole is a 
lorentzian manifold $(M, g)$ solution of the Einstein 
equation\footnote{The full Einstein equation contains 
an additionnal term on the right side, which is a symmetric 
trace-free tensor describing the matter and physical forces 
in the spacetime. The equation stated here is an equation in the void.}:

\[  \mbox{Ric}_g - \frac{R}{2}g =\Lambda.g \]

where:

-- $\mbox{Ric}_g$ is the Ricci tensor,

-- $R$ is the scalar curvature, i.e. the trace of the Ricci tensor 
(relatively to the metric $g$),

-- $\Lambda$ is a prescribed  real number, the \emph{cosmological constant.\/}

It was soon realized by K. Schwarzschild that this equation 
with $\Lambda= 0$ admits the following $1$-parameter family of solutions:

\[ ds^2 = - (1 - \frac{2M}{r}) dt^2 + \frac{dr^2}{1-2M/r} + r^2ds_{0}^{2} \]

In this expression, $ds_0^2 = d\theta^2+\sin^2\theta d\phi^2$ is the usual 
round metric on the sphere ${\mathbb S}^2$. The term $r$ is 
$\sqrt{x^2+y^2+z^2}$, where $(x,y,z,t)$ are coordinates of 
${\mathbb R}^4$, and $(\theta, \phi)$ are the spherical coordinates 
of $(\frac{x}{r}, \frac{y}{r}, \frac{z}{r}) \in {\mathbb S}^2$.

Moreover, in some manner, the Schwarzschild metric is the 
``unique'' solution of the Einstein equation on ${\mathbb R}^4$ with 
spherical symmetry, i.e. invariant by the usual 
$\mbox{SO}(3,{\mathbb R})$-action (Birkhoff Theorem, see 
\cite{gravitation} for a more rigorous and detailed exposition). 
Hence, this metric suits perfectly the role of model of a 
non-rotating isolated object (for example, a star) with mass 
$M$. But these solutions contains a psychological difficulty: 
the radius of the object could be inferior to $r=2M$, in which case the
singularity $r=2M$ of the metric cannot be dropped out by considering 
that it is hidden inside the object, where the general relativity 
theory does not apply anymore. It has been realized that the 
singularity is not a real one: it appears only in reason of the 
selected coordinate system. Actually, the Schwarzschild metric 
on ${\mathbb R}^4 \setminus \{ r = 0, 2M \}$ can be isometrically 
embedded in another lorentzian manifold $M_{KS}$, diffeomorphic to 
${\mathbb R}^2 \times {\mathbb S}^2$, and maximal: $M_{KS}$ cannot 
be isometrically embedded in a bigger spacetime.  This maximal 
spacetime is usually described in the following way: equip 
${\mathbb R}^2 \times {\mathbb S}^2$ with coordinates 
$(u, v, \theta, \phi)$ - the \emph{Kruskal-Szekeres coordinates\/} - 
where $(\theta, \phi)$ are still the spherical coodinates.
Then, $M_{KS}$ is the open domain $U \times {\mathbb R}^2$, where 
$U = \{ v^2 - u^2 < 1 \}$, equipped with the metric:

\[ ds^2_{KS} = (32M^3/r)e^{-r/2M}(-dv^2+du^2) + r^2 ds^2_{0} \]

The term $r$, defined by $(r/2M - 1)e^{r/2M} = u^2 -v^2$ coincide 
with the $r$-coordinate on 
the Schwarzschild domain: hence, the ``singularity'' $r=2M$ vanishes: 
it corresponds to the locus $\{ u^2 = v^2 \}$.

The domain $O = \{ 0 < u^2 - v^2 <1\}$ is the \emph{outer domain:\/} it 
is thought as the region where the typical (prudent) observer takes place. 
There are two other regions: $B^{\pm} = \{ u^2 - v^2 < 0, \pm v > 0 \}$. 
Moreover, $M_{KS}$ is time-oriented so that $v$ increases with time.

$B^+$ enjoys the following remarkable property: there is no future-oriented 
causal curve starting from a point in $B^+$ and reaching a point 
in $O$. Similarly, a past-oriented causal curve cannot escape 
from $B^-$. In a more ``physical'' language, photons cannot escape 
from $B^+$: $B^+$ is invisible from $O$. In other words, $B^+$ is a 
\emph{black-hole.\/} The no-return frontier $\{ v = \mid{u}\mid \}$ 
is the \emph{(event) horizon.\/}

The family of Schwarzschild metrics is actually included in a more general 
family of solutions, the \emph{Kerr metrics:\/}

\[ -(1-\frac{2Mr}{\rho^2})dt^2 + \frac{\rho^2}{\Delta}dr^2 - 
{\sin}^2(\theta)\frac{4Mra}{\rho^2}d\phi{dt} + \rho^2d\theta^2 + 
(r^2+a^2+ {\sin}^2(\theta)\frac{2Mra}{\rho^2})\sin^2(\theta) d\phi^2 \]

where:

-- $\Delta = r^2 - 2Mr +a^2$, t $\rho^2 = r^2 + a^2\cos^2(\theta)$,

-- $M$ is positive (the ``mass''),

-- $a$ is a real number (the ``angular momenta per mass unit''). 

A phenomenom similar to the embedding of Schwarzschild spacetime in 
Kruskal-Szeke\-res coordinates still applies:  every Kerr-spacetime can be embedded 
in a na\-tu\-ral way into a maximal spacetime $M_{Kerr}^{max}$ where the singularities 
appearing in the Kerr-coordinates vanish, and where some domains deserve the 
appelation ``black-hole''. But the description of this completion is 
much more involved: it is the entire matter of the excellent book \cite{Oneillkerr}, 
which is also one of the best references for a rigorous description of Schwarzschild 
black-holes. 




\subsection{Towards a general definition of black-holes}
The basic facts on Schwarz\-schild spacetimes presented above are 
sufficient to provide a quite satisfactory illustration of essential 
notions that should appear in any ma\-thematical definition of black-holes.
Let's be more precise and express some problems arising when trying 
to elaborate such a general theory, i.e. in our meaning, to specify 
what are the spacetimes  deserving the appelation of spacetime with black-hole. 

\begin{enumerate}

\item \emph{Where are the typical observers?\/} 
The description of the Schwarzschild black-hole makes clear that this 
notion is relative to the region $O$ where the observers are assumed to stay. 
A black-hole is simply a connected component of the region invisible from $O$.
How to define the region $O$ in $M$ without specifying the black-hole itself?

\item \emph{Is some part of $M$ missing?\/} The notion of black-hole is a global 
property, depending on the entire $M$. For example, if we delete from the Minkowski 
space some regions,
we can easily produce in an artificial way regions invisible from the observers 
(assuming solved the first question above) which does not correspond to a 
physically relevant example of black-hole. Hence, we have to define a notion 
of ``full'' spacetime ensuring that some part of the invisible domain is not 
simply due to the absence of some relevant region.
\end{enumerate}

\rque
\label{rk.censure}
To these two basic problems, another important requirement, 
traditionnally appearing in the physical litterature, should be added: 
the \emph{Cosmic Censorship.\/} This condition admits many different formulations. 
It is most of the time 
expressed in the form ``singularities of the spacetime must be hidden to 
the observers by the horizons''. There is a geometrical way to translate 
this notion, interpreting singularities as final points of non-complete 
causal geodesics, or as causal singularities: we could define 
spacetimes satisfying the Cosmic Censorship if causal closed 
curves, if any, all belong to a black hole, and if any non-complete 
causal future-oriented ray in $M$, with starting point visible 
from one observer, must enter in a black-hole

There is a particular case of the  Cosmic Censorship conjecture 
often stated as follows: ``generically'' (to be defined!), maximal 
globally hyperbolic spacetimes are maximal, i.e. cannot be embedded 
in a bigger spacetime, even not globally hyperbolic. 
But the spacetimes containing a black-hole we will consider here
are not globally hyperbolic, 
hence this expression of Cosmic Censorship is not relevant for them. 
\erque

Some answer to the first question above seems to be 
widely accepted: 
\emph{observers lie near the conformal boundary.\/}  
Then, we can define the black-holes as the connected components of the 
\emph{invisible domain,\/} i.e. the interior of the region of $M$ 
containing the points $x$ such that no causal future-oriented 
curve starting from $x$ tends to a point in the conformal boundary of $M$. 
Of course, this solution arises immediately another question: what is the 
conformal boundary? This new question is not answered in full generality - for
a nice recent survey on related questions see \cite{senovilla}.
We can summarize many works by observing that this point is most of the 
time solved simply by prescribing from the beginning what is the conformal 
boundary, without making sure that the proposed boundary is maximal in any 
meaning\footnote{Anyway, for the Kruskal-Szekeres spacetimes $M_{KS}$ and 
maximal Kerr-spacetimes $M_{Kerr}^{max}$, which are all asymptotically flat, 
the conformal boundary is well-defined in a fully satisfying way.}.
However, in our special AdS context, 
we have a completely satisfactory and easy definition of conformal boundary 
of \emph{strongly causal\/} AdS spacetimes (see next {\S}).

Concerning the second question, encouraged by the examples of Kruskal-Szekeres 
and Kerr spacetimes, one could hope that a good answer is simply to require $M$ 
to be maximal in the sense that it cannot be isometrically embedded in a bigger 
spacetime satisfying the Einstein equation. Unfortunately, this attempt fails 
for BTZ black-holes: they are not maximal in this meaning. Furthermore, such 
a definition of full would lead to some incoherence in the maximal Kerr spacetime.
The restriction to causal spacetimes, i.e. to require that 
$M$ is a causal spacetime, and that it cannot be embedded in a bigger causal 
spacetime does not solve the problem: once more, it does not 
apply to BTZ black-holes (see remark~\ref{rk.pasmaximal}). 

Our deceiving conclusion is that we still don't know how to express 
in a satisfying way what is a good notion of full spacetime, even 
in the AdS-background. As a positive element of answer, we stress 
out that all the examples reproduced later enjoy the following 
properties of ``maximal'' nature: every connected component of 
their conformal boundary is maximal globally hyperbolic,
and every black-hole, i.e. every connected component of the 
invisible domain (but this domain will always be connected) 
is maximal globally hyperbolic too (see remark~\ref{rk.bordcausalgh}).
Observe that these 
requirements (global hyperbolicity of black-holes and 
observer-spacetime) reflect in some way the Cosmic 
Censorship principle for globally hyperbolic spacetimes 
(remark~\ref{rk.censure}).

\subsection{BTZ black-holes}
First of all, BTZ-multi-black-holes have dimension $2+1$, i.e. 
are $3$-dimensional. In this low dimension the Einstein equation 
is remarkably simplified: the solutions have all constant sectional 
curvature with the same sign that the cosmological constant $\Lambda$. 
We only consider the case $\Lambda < 0$. Hence, up to rescaling, 
the spacetimes satisfying the Einstein equation are precisely the spacetimes 
locally modeled on AdS.

This special feature allows us to propose a correct answer to the first question 
in the preceding paragraph, in the $2+1$-dimensional case, i.e. with 
AdS-background, thanks to the natural conformal completion of AdS 
by $\mbox{Ein}_2$. More precisely, the spacetimes $M$ we will consider 
are locally modeled in AdS - in short, they are AdS-spacetimes. Hence, 
their universal covering admits a development 
${\mathcal D}: \widetilde{M} \rightarrow \mbox{AdS}$. Define then the lifted 
conformal boundary $\Omega$ as the interior of the set of final extremities 
in $\mbox{Ein}_2 =\partial{\mbox{AdS}}$ of future-oriented causal curves 
with relative interior contained in ${\mathcal D}(\widetilde{M})$. We could 
try to define the conformal boundary of $M$ as the quotient of 
$\Omega$ by the holonomy group $\rho(\Gamma)$. Unfortunately, this definition 
admits an uncomfortable drawback: the development ${\mathcal D}$ could be 
non-injective, which requires a modification of the definition above. 
Moreover, we have no guarantee \emph{a priori\/}  that the action of 
$\rho(\Gamma)$ on $\Omega$ is free and proper. 

In \S~\ref{sub.bord}, we show to solve these difficulties in the context of
strongly causal spacetimes.: every 
strongly causal AdS-spacetime admits a natural strongly causal conformal boundary 
(Theorem~\ref{thm.bordcausal}).

For that reason (besides the 
physical coherence of such an assumption) we restrict from now to 
strongly causal spacetimes. Now we can state our peculiar definition:

\begin{defin}
\label{def.BH}
An AdS-spacetime with black hole is a strongly causal AdS spacetime $M$ such that:

-- $M$ admits a non-empty strongly causal conformal boundary $O$,

-- the past of $O$, i.e. the region of $M$ made of initial points of future oriented
causal curves ending in $O$,is not the entire $M$.

Every connected component of the interior of the complement in $M$ of the past of $O$
is a black-hole.
\end{defin}

Let's now collect the examples. While considering discrete groups of isometries of AdS, 
we have constructed many examples of spacetimes locally modeled on AdS: the quotients 
$M(\Gamma) = \Gamma\backslash{E}(\Gamma)$ where $\Gamma$ is a non-abelian torsion-free 
discrete admissible subgroup of Isom(AdS).  \emph{Every $M(\Gamma)$ contains a black-hole, 
except if it is globally hyperbolic, i.e. except if  $\Lambda$ is a topological circle.\/} 

This construction applies in more cases, even if $\Gamma$ is abelian, 
but the elementary achronal subset $\Lambda$ must be selected: it must contain at least 
two points, and $\Gamma$ must be a cyclic subgroup, generated by an isometry $\gamma_0$. 
More precisely:

-- (the conical case) This is the case where $\Lambda$ is 
the union of two lightlike geodesic segments $[y,x]$, $[z,x]$, 
each of them not reduced to single point and with a common 
extremity $x$, which is their common past. Then
$\gamma_0$ must be a hyperbolic - hyperbolic element 
(see definition~\ref{def.synch}) and $y$, 
$z$ must be attractive fixed points of $\gamma_0$.

-- (the splitting case) This is the case where $\Lambda$ 
is a pair of non-causally related points in $\mbox{Ein}_2$. 
As in the previous case, these points must be the attractive 
and repulsive fixed points of the hyperbolic-hyperbolic element
$\gamma_0$.

-- (the extreme case) This last case is the case where $\Lambda$ 
is a lightlike segment, not reduced to a single point. Then, 
$\gamma_0$ must be a parabolic-hyperbolic element 
(see definition~\ref{def.synch}) fixing the two extremities of $\Lambda$. 

In all these cases, the quotient $M_\Lambda(\Gamma) = \Gamma\backslash{E}(\Lambda)$ 
is still a spacetime with conformal boundary, and with non-empty invisible set: a black-hole.

Since it requires a basic knowledge of AdS geometry, the description of these spacetimes 
is postponed to {\S}~\ref{sec.last}. We just mention in this introduction that they 
essentially\footnote{Except the particular the case of non-static single BTZ black-holes, 
see remark~\ref{rk.singleBTZ}.} include all the previously examples named BTZ black-holes 
(\cite{BTZ, BTZ2, amaking, aminneborg, amispin, brill1, brill}, etc...), but contains also 
new examples: mainly, the case where  $\Gamma$ is not finitely generated, and also the 
conical case, which is not considered as a spacetime with black-hole in these references, 
probably because it is obviously non-maximal, since, in this case, $M_\Lambda(\Gamma)$ 
embeds isometrically into $M_{yz}(\Gamma)$.

\section*{Acknowledgements}
L. Freidel introduced me to the realm of BTZ black-holes and 
their physical motivations. He initiated my special interest 
to Lorentzian geometry with constant curvature. F. B\'eguin, 
C. Frances, A. Zeghib have participated to the early stage of 
the conception of this work in the mathematical and geometrical 
understanding of BTZ-black holes and related questions. 
P. Saad\'e has performed the pictures.

\section*{Notations}
This paper should be readen jointly with \cite{ba1} in which
all basic notions and notations are introduced. Let's recall
some of them, for the reader convenience:

\begin{itemize}

\item AdS is the anti-de Sitter space, i.e. the locus $\{ Q = -1 \}$
in $E \approx {\mathbb R}^{2,2}$ equipped with a quadratic form $Q$ of
signature $(2,2)$. It is isometric to
the Lie group $\mbox{SL}(2, {\mathbb R})$ endowed with the lorentzian
metric defined by the Killing form. Its universal covering
is denoted by $\widetilde{\mbox{AdS}}$ and its projection
into the sphere $S(E)$ of rays in $E$ - the Klein model - is denoted 
${\mathbb A}{\mathbb D}{\mathbb S}$. The projective Klein model
$\overline{{\mathbb A}{\mathbb D}{\mathbb S}}$ is the projection 
in the projective space $P(E)$. It is naturally identified with
the Lie group $G = \mbox{PSL}(2, {\mathbb R})$.
All these lorentzian manifolds are oriented and chronologically oriented.

\item A spacetime is a lorentzian manifold which is oriented and chronologically
oriented. An AdS - spacetime is a spacetime locally modelled on the anti-de Sitter
space AdS.

\item Affine domains in AdS, $\widetilde{\mbox{AdS}}$, 
${\mathbb A}{\mathbb D}{\mathbb S}$ are lifts of affine domains in 
$\overline{{\mathbb A}{\mathbb D}{\mathbb S}}$, i.e. intersection 
in $P(E)$ between $\overline{{\mathbb A}{\mathbb D}{\mathbb S}}$ and
affine domains $U$ of $P(E)$ such that the intersection between the 
projective plane $\partial{U}$ and 
$\overline{{\mathbb A}{\mathbb D}{\mathbb S}}$ is a spacelike surface; it is then
an isometric copy of the hyperbolic plane ${\mathbb H}^2$. 

\item The isometry group of $\overline{{\mathbb A}{\mathbb D}{\mathbb S}}$ is 
$G \times G$. The isometry group of $\widetilde{\mbox{AdS}}$ is isomorphic
to $\widetilde{G} \times \widetilde{G}$ quotiented by the cyclic diagonal
group generated by $(\delta, \delta)$, where $\delta$ is a generator
of the center of $\widetilde{G}$, the universal covering of $G$.

\item There is a conformal embedding of AdS into the Einstein universe 
$\mbox{Ein}_3$, which is ${\mathbb S}^2 \times {\mathbb S}^1$ equipped with
(the conformal class of) the lorentzian metric $ds^2 - dt^2$, where
$ds^2$ is the round metric on the $2$-sphere and $dt^2$ the usual metric
on ${\mathbb S}^1 \approx {\mathbb R}/{2\pi{\mathbb Z}}$.
The universal covering of $\mbox{Ein}_3$ is denoted by $\widehat{\mbox{Ein}}_3$.
There is a conformal embeding of $\widetilde{\mbox{AdS}}$ into 
$\widehat{\mbox{Ein}}_3$. The boundary of the image of this embedding is the conformal 
boundary $\partial\widetilde{\mbox{AdS}} \approx \widehat{\mbox{Ein}}_2$.
The projection $\overline{\mbox{Ein}}_2$ of $\mbox{Ein}_2$ in $P(E)$ 
is the conformal boundary of $\overline{{\mathbb A}{\mathbb D}{\mathbb S}}$.

\item $\overline{\mbox{Ein}}_2$ is bifoliated by two transverse foliations 
$\overline{\mathcal G}_L$, $\overline{\mathcal G}_R$. 
Every leaf of $\overline{\mathcal G}_L$ (the left foliation) or of 
$\overline{\mathcal G}_R$ (the right foliation) is a lightlike
geodesic, canonically isomorphic to the real projective line ${\mathbb R}P^1$. 
Every leaf of $\overline{\mathcal G}_L$ intersects every leaf of 
$\overline{\mathcal G}_R$ at one and only one point. Finally, the 
leaf space of the left (resp. right) foliation is naturally isomorphic
to the real projective line. It is denoted by 
${\mathbb R}P^1_L$ (resp. ${\mathbb R}P^1_R$).

\item For all causality notions, we refer to \cite{beem}. In \cite{ba1}
we also present the causality notions pertinent for our purpose. An achronal
subset of $\widehat{\mbox{Ein}}_2$ is a subset $\widetilde{\Lambda}$
such that for any pair $(x,y)$ of distinct elements of $\widetilde{\Lambda}$,
$x$ and $y$ are not causally related in $\widehat{\mbox{Ein}}_2$. An achronal
subset $\widetilde{\Lambda}$ is pure lightlike if it contains two opposite
elements $x$, $\delta(x)$; if not, it is generic. A generic achronal subset
is elementary if it is contained in an union $l \cup l'$, where $l$, $l'$
are lightlike geodesic segments in $\widehat{\mbox{Ein}}_2$ (see 
{\S} $8.7$ in \cite{ba1}).
All these notions project on similar notions in $\mbox{Ein}_2$. See \cite{ba1},
{\S} $5$ for more details.

\item The invisible domain of a generic closed achronal subset 
$\widetilde{\Lambda}$ in $\widehat{\mbox{Ein}}_2$ is the set of points
in $\widehat{\mbox{Ein}}_2$ which are not causally related to any element 
of $\widetilde{\Lambda}$. It is denoted by $\Omega(\widetilde{\Lambda})$.
It is the region between the graphs $\widetilde{\Lambda}^+$ and 
$\widetilde{\Lambda}^-$ of two $1$-Lipschitz functions from ${\mathbb S}^1$ into
$\mathbb R$ (see {\S}~$8.6$ in \cite{ba1}).
The invisible domain $E(\widetilde{\Lambda})$ in $\widetilde{\mbox{AdS}}$ is
the set of points in $\widetilde{\mbox{AdS}} \subset \widehat{\mbox{Ein}}_3$
which are not causally related to any element 
of $\widetilde{\Lambda}$ in $\widehat{\mbox{Ein}}_3$.

\end{itemize}

\section{Strongly causal spacetimes}
\label{sec.strong}

A lorentzian manifold $M$ is \emph{strongly causal\/} if 
for every point $x$ in $M$ every neighborhood of $x$ 
contains an open neighborhood $U$ which 
is \emph{causally convex,\/} i.e. such that any causal curve in $M$ 
joining two points in $U$ is actually contained in $U$
(see \S $2.3$ in \cite{ba1}). An
isometric action of a group $\Gamma$ on a strongly causal lorentzian
manifold $M$ is strongly causal if any point $x$ of $M$ admits 
arbitrarly small neighborhoods $U$ such that for every non trivial element $\gamma$ of
$\Gamma$ no point of $U$ is causally related to a point of $\gamma{U}$.
Clearly, if $\Gamma$ is a group of isometries of a strongly causal lorentzian
manifold $M$ acting freely and properly discontinuously, the quotient
manifold $\Gamma\backslash{M}$ is strongly causal if and only if
the action of $\Gamma$ is strongly causal.

In this paper we need to prove that certain isometric actions are strongly
causal. In this {\S} we provide two lemmas useful for this purpose.

In the first lemma we consider a strongly causal lorentzian manifold $M$
admitting a Killing vector field
$X$. We assume that $X$ is everywhere spacelike. Let $\phi^t$ be the flow generated by
$X$ and assume that the ${\mathbb R}$-action defined by $\phi^t$ is free and
properly discontinuous. Let $Q^\phi$ be the orbit space of $\phi^t$ and
$\pi: M \rightarrow Q^\phi$ the projection map. 
It is easy to show that $\pi$ is a locally trivial fibration with fibers
homeomorphic to $\mathbb R$.
For every tangent vector $v$
of $Q^\phi$ at a point $y$ let $x$ be any element of $\pi^{-1}(y)$ and let $w$ be the
unique tangent vector to $M$ at $x$ orthogonal to $X(x)$ and such that $d\pi_x(w) = v$. 
Since $\phi^t$ is isometric the norm of $w$ does not depend on the choice of $x$:
we define the norm $\mid v \mid$ as this norm $\mid w \mid$. This procedure defines
a lorentzian metric on $Q^\phi$.

\begin{lem}
\label{lem.quocausal}
If the lorentzian manifold $Q^\phi$ is strongly causal then the action of $\phi^1$ on $M$ is
strongly causal.
\end{lem}

\preu
Let $x$ be any element of $M$ and $U$ any neighborhood of $x$.
Let $y = \pi(x)$. Since $Q^\phi$ is assumed strongly
causal, $U' = \pi(U)$ contains a causally convex open
neighborhood $V'$ of $y$. Define $V = \pi^{-1}(V')$: this open subset is
diffeomorphic to $V' \times {\mathbb R}$, such that the restriction of 
$\phi^t$ to $V$ corresponds to the translation on the $\mathbb R$-factor. 
We also take the convention that the $\mathbb R$-coordinate of $x$ is
zero: $x \approx (y,0)$. Shrinking $V'$ if necessary, we can assume that
the spacelike norm of $X$ on $V$ is uniformly bounded by some constant $C$,
and moreover that any causal curve in $V'$ has lorentzian length less than
$1/2C$. It then follows that the $\mathbb R$-coordinate of any causal
curve in $V$ cannot vary more than $1/2$. Let $W$ be a causally convex
open neighborhood of $x$ contained in
$V' \times ]-1/4, +1/4[ \subset V$.
We claim that for every non zero integer $n$ non point of $\phi^nW$ is
causally related to a point of $W$. Indeed: assume by contradiction
that there is a causal curve $c: I \rightarrow M$ joining an
element $(z,t)$ to an element $(z', s+n)$ with $z$, $z'$ in $V'$ and
$s$, $t$ in $]-1/4, +1/4[$. The projection $\pi \circ c$ is a causal
curve in $Q$ joining two $\pi(z)$, $\pi(z')$. Since $V'$ is causally convex
$\pi \circ c$ is contained in $V'$, i.e. $c$ is contained in $V$. Hence
the $\mathbb R$-coordinate cannot vary along $c$ more than $1/2$.  
We obtain a contradiction since $s+n-t$ has absolute value bigger than $1/2$.\fin

\rque
\label{rk.quocausalfini}
Lemma~\ref{lem.quocausal} remains true when the action of $\phi^t$ is not free
but periodic. In this case the group generated by $\phi^1$ is finite.
Details are left to the reader.
\erque

We will apply lemma~\ref{lem.quocausal} in the case where $M$ has dimension $3$ and is simply
connected. In this situation the map $\pi$ is a fibration and the homotopy sequence of
this fibration implies that $Q^\phi$ is a $2$-dimensional manifold which is simply connected.
Hence a very nice complement is:

\begin{lem}
\label{lem.2causal}
Any simply connected $2$-dimensional lorentzian manifold is strongly causal.
\end{lem}

\preu
Any simply connected surface is diffeomorphic to the sphere or the plane ${\mathbb R}^2$.
Since the sphere does not admit any lorentzian metric we just have to consider the
case of the plane ${\mathbb R}^2$. Any reader acquainted with dynamical systems
will recognize that the lemma follows quite easily from
the Poincar\'e-Bendixon Theorem applied to the lightlike foliations: 
every leaf of one of these foliations is a closed embedding of the real line
in ${\mathbb R}^2$ and for any point $x$ in ${\mathbb R}^2$ the future (resp. 
the past) of $x$ is the domain bounded by $r_1 \cup r_2$ where $r_1$, $r_2$ are
the future oriented (resp. past oriented) lighlike geodesic rays starting from $x$.
Details are left to these readers.

For other readers, we refer to Theorem $3.43$ of \cite{beem} where 
a slightly better statement is proved:
lorentzian metrics on the $2$-plane are \emph{stably\/} causal, which
is stronger than being strongly causal.
\fin

\section{Description of the elementary generic achronal subsets of 
$\widehat{\mbox{Ein}}_2$}
\label{sub.elementary}

In this section, we complete the descriptions of invisible domains 
$E(\widetilde{\Lambda})$ and $\Omega(\widetilde{\Lambda})$ by considering
the elementary cases. For the reader convenience, we start with the splitting 
case, which is already described in \cite{ba1}, {\S} $8.8$.

\subsection{The splitting case}
\label{subsub.split}
It is the case $\widetilde{\Lambda} = \{ x,y \}$, where $x$, $y$ are 
two non-causally related points in $\widehat{\mbox{Ein}}_{2}$. 
Then, $\{ x,y \}$ is a gap pair and there are two associated ordered 
gap pairs that we denote respectively by $(x,y)$ and $(y,x)$ (see definition
$8.24$ in \cite{ba1}). 
$\widetilde{\Lambda}^+$ is an union 
${\mathcal T}^+_{xy} \cup {\mathcal T}^+_{yx}$ of two future oriented lightlike 
segments with extremities $x$, $y$ that we call \emph{upper tents\/.} 
Such an upper tent is the union of two lightlike segments, one starting 
from $x$, the other from $y$, and stopping at their first intersection point, 
that we call the \emph{upper corner.\/}

\begin{figure}[ht]
\centerline{\includegraphics[width=8cm]{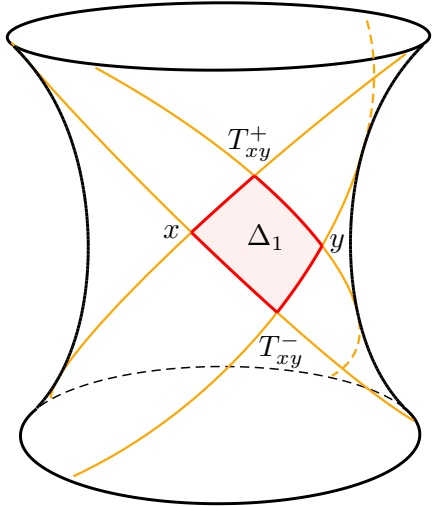}}
\caption{\label{ttents}
   \textit{Upper and lower tents}} 
\end{figure}

\begin{figure}
\centerline{\includegraphics[width=8cm]{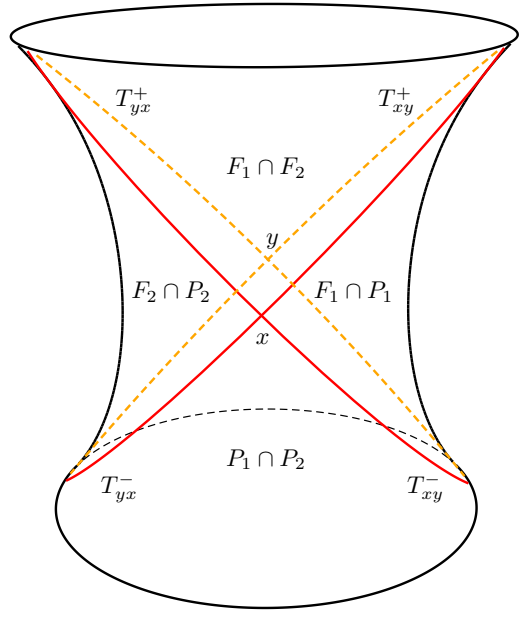}}
\caption{\label{splitfigure}
   \textit{The splitting case. The domain $E(\Lambda)$ is between 
the hyperplanes $x^\perp$ and $y^\perp$. These hyperplanes,
tangent to the hyperboloid, are not drawn, except their intersections with the hyperboloid,
which are the upper and lower tents ${\mathcal T}^\pm_{xy}$, ${\mathcal T}_{yx}^\pm$.}} 
\end{figure}

Similarly, $\widetilde{\Lambda}^-$ is an union 
${\mathcal T}^-_{xy} \cup {\mathcal T}^-_{yx}$ of two 
\emph{lower tents\/} admitting a similar 
description, but where the lightlike segments starting from $x$, $y$ 
are now past oriented
(see Figure \ref{ttents}, Figure~\ref{splitfigure}), and sharing a common extremity: the \emph{lower 
corner.\/}

The invisible domain $\Omega(\widetilde{\Lambda})$ from $\widetilde{\Lambda}$ in $\widehat{\mbox{Ein}}_{2}$ is the union of two diamond-shape regions $\widetilde{\Delta}_1$, $\widetilde{\Delta}_2$. The boundary of $\widetilde{\Delta}_1$ is the union ${\mathcal T}^+_{xy} \cup {\mathcal T}^-_{xy}$, and the boundary of $\widetilde{\Delta}_2$ is ${\mathcal T}^+_{yx} \cup {\mathcal T}^-_{yx}$. We project all the picture in some affine region $V \approx {\mathbb  R}^3$ of $S(E)$ such that:

-- $V \cap {\mathbb A}{\mathbb D}{\mathbb S}$ is the interior of the hyperboloid: $\{ x^2 + y^2 < 1 + z^2 \}$,

-- $\Lambda = \{ (1,0,0), (-1, 0, 0) \}$.

Then, $E(\Lambda)$ is region $\{ -1 < x < 1 \} \cap {\mathbb A}{\mathbb D}{\mathbb S}$.
One of the diamond-shape region $\widetilde{\Delta}_i$ projects to $\Delta_1 = \{ -1 < x < 1, y >0 , x^2+y^2 = 1+z^2 \}$, the other projects to $\Delta_2 = \{ -1 < x < 1, y <0 , x^2+y^2 = 1+z^2 \}$. The past of $\Delta_1$ in $E(\Lambda)$ is $P_1 = \{ (x,y,z) \in E(\Lambda) / z < y \}$. and the future of $\Delta_1$ in $E(\Lambda)$ is $F_1 = \{ (x,y,z) \in E(\Lambda) / z > -y \}$. We have of course a similar description for the future $F_2$ and the past $P_2$ of $\Delta_2$ in $E(\Lambda)$. Observe:

-- the intersections $F_1 \cap F_2$ and $P_1 \cap P_2$ are disjoint. They are tetraedra in $S(E)$:
$F_1 \cap F_2$ is the interior of the convex hull of $\Lambda^+$, 
and $P_1 \cap P_2$ is the interior of the 
convex hull of $\Lambda^-$.

-- the intersection $F_1 \cap P_1$ (resp. $F_2 \cap P_2$) is the intersection between ${\mathbb A}{\mathbb D}{\mathbb S}$ and the interior of a tetraedron in $S(E)$: the convex hull of $\Delta_1$ (resp. $\Delta_2$).

\begin{defin}
\label{defghsimple}
$E^+(\Lambda) = F_1 \cap F_2$ is the future globally hyperbolic convex core; $E^-(\Lambda) = P_1 \cap P_2$ is the past globally convex core.
\end{defin}

This terminology is justified by the following (easy) fact: $F_1 \cap F_2$ (resp. $P_1 \cap P_2$) is the invisible domain $E(\Lambda^+)$ (resp. $E(\Lambda^-)$). Hence, they are indeed globally hyperbolic.

The intersection between the closure of $E(\Lambda)$ in $S(E)$ and the boundary ${\mathcal Q}$ of ${\mathbb A}{\mathbb D}{\mathbb S}$ is the union of the closures of the diamond-shape regions. Hence, $\Delta_{1,2}$ can be thought as the conformal boundaries at infinity of $E(\Lambda)$. Starting from any point in $E(\Lambda)$, to $\Delta_i$ we have to enter in $F_i \cap P_i$, hence we can adopt the following definition:

\begin{defin}
\label{def.end}
$F_1 \cap P_1$ is an end of $E(\Lambda)$.
\end{defin}

Finally:

\begin{defin}
The future horizon is the past boundary of $F_1 \cap F_2$; the past horizon is the
future boundary of $P_1 \cap P_2$.
\end{defin}

\begin{prop}
\label{simpledecomp}
$E(\Lambda)$ is the disjoint union of the future and past 
glo\-bally hyperbolic cores $E^\pm(\Lambda)$, of the two ends,  
and of the past and future horizons.\fin
\end{prop}

\rque
\label{rk.BTZnotations}
In the conventions of \cite{BTZ, BTZ2, brill}, the globally hyperbolic convex cores $F_1 \cap F_2$ and $P_1 \cap P_2$ are \emph{regions of type II,\/} also called \emph{intermediate regions.\/} The ends $F_1 \cap P_1$ and $F_2 \cap P_2$ are \emph{outer regions,\/} or \emph{regions of type I.\/}
\erque

\subsection{The extreme case}
\label{subsub.extreme}

\begin{figure}
\centerline{\includegraphics[width=8cm]{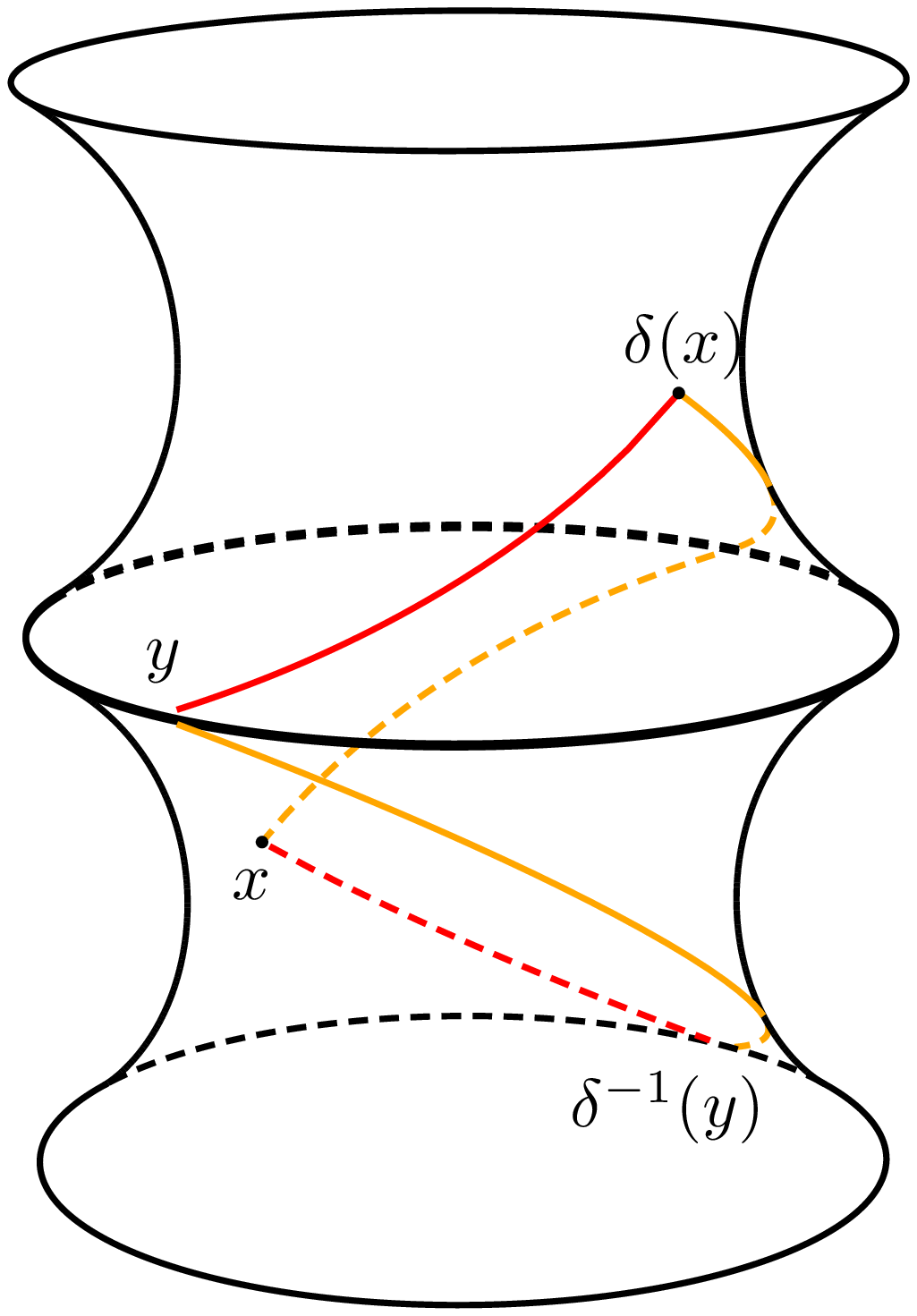}}
\caption{\label{extremefigure}
   \textit{The extreme case. The point $\delta(x)$ is in the front of
   the upper hyperboloid. The point $x$ is in the rear of the lower hyperboloid.
   The domain $E(\Lambda)$ is the domain between the hyperplanes $\delta(x)^\perp$ and 
   $\delta^{-1}(y)^\perp$.}} 
\end{figure}

The extreme case is harder to picture out 
since $\Omega(\widetilde{\Lambda})$ and $E(\widetilde{\Lambda})$ 
are \emph{not} contained in an affine domain (see figure~\ref{extremefigure}). 
Assume that $y$ is in the future of $x$. Observe that $\widetilde{\Lambda}^\pm$ 
are then pure lightlike. Hence, $E(\widetilde{\Lambda}^\pm)$ are empty. 
The region $\Omega(x,y)$ is a ``diamond'' in $\mbox{Ein}_2$ 
(we call it an \emph{extreme diamond\/} ) admitting as 
boundary four lightlike segments: the segments $[y, \delta(x)]$, 
$[x, \delta(x)]$, $[\delta^{-1}(y), x]$, and $[y, \delta^{-1}(y)]$. 

A carefull analysis shows that $E(x,y)$ is precisely the intersection 
between the past and the future of $\Omega(x,y)$.

Keeping $x$ fixed, and considering a sequence $y_n$ converging to $y$, 
with $y_n$ non-causally related to $x$, then one of the diamond-shape 
region $\Delta^n_i$ of the associated $\Omega(x, y_n)$ - let's say, 
$\Delta^n_2$ - vanishes. The other converges to the entire region 
$\Omega(x,y)$. The various parts of the domains $E(x, y_n)$, namely the 
globally hyperbolic convex cores $E^\pm(x,y_n)$ and the ends, vanish, 
except one of the end, which becomes closer and closer to the entire $E(x,y)$. 

\subsection{The conical case}
\label{subsub.conical}
In the conical case, $\Lambda$ is an upper or lower tent. By symetry, 
we can consider only the upper case: 
$\Lambda  ={\mathcal T}^+_{xy} = [x,z] \cup [z,y]$. Then, 
$\Lambda^+ = \Lambda \cup {\mathcal T}^+_{yx}$, and 
$\Lambda^-$ is the pure lightlike subset 
$\Lambda \cup {\mathcal T}^-_{yx}$. In the notations of \S~\ref{subsub.split}, 
$\Lambda^+$ is the future boundary of the diamond $\Delta_1$, and 
$\Omega(\Lambda)$ is the diamond $\Delta_2$.
Then, in some affine domain, $E(\Lambda) = \{ (x,y,z) / -1 < x < 1, z>y \}$. 
It can be described as the future in $V$ of $\Delta_2$. 
It is also the intersection between the past of $x$, 
the past of $y$, and the complement of the past of $z$. Finally, 
$E(\Lambda)$ is the union of $F_1 \cap F_2$, $F_2 \cap P_2$ and 
the component of the past horizon of $E(x,y)$ separating these two regions 
(see figure \ref{splitfigure}).

\section{Causal domains of isometries of AdS  }
\label{sec.cyclic}

\subsection{The isometry group}
We recall some facts established in \cite{ba1}, {\S} $9$, concerning isometries 
preserving generic
achronal subsets of $\mbox{Ein}_2$.
We use the identification 
$\overline{{\mathbb A}{\mathbb D}{\mathbb S}} \approx G = \mbox{PSL}(2, {\mathbb R})$ 
(cf. {\S} Notations). Then $\widetilde{\mbox{AdS}}$ can be identified with the universal covering 
$\widetilde{G} = \widetilde{\mbox{SL}}(2, {\mathbb R})$. Denote by 
$\bar{p}: \widetilde{G} \rightarrow G$ 
the covering map, and $Z$ the kernel of $\bar{p}$: $Z$ is cyclic, it is the center of 
$\widetilde{G}$. 
Let $\delta$ be a generator of $Z$: we select it in the future of the neutral element $id$.

$\widetilde{G} \times \widetilde{G}$ acts by left and right translations on $\widetilde{G}$.
This action is not faithfull: the elements acting trivially are precisely the elements in 
$\mathcal Z$, the image of $Z$ by the diagonal embedding. The isometry group 
$\widehat{\mbox{SO}}_0(2,2)$ is then identified with 
$(\widetilde{G} \times \widetilde{G})_{/\mathcal Z}$.

Let $\mathcal G$ be the Lie algebra $\mbox{sl}(2, {\mathbb R})$ of $G$: the Lie algebra
of $(\widetilde{G} \times \widetilde{G})_{/\mathcal Z}$ is ${\mathcal G} \times {\mathcal G}$. 
We assume the reader familiar with the notion of elliptic, parabolic, hyperbolic elements of 
$\mbox{PSL}(2, {\mathbb R})$. Observe that hyperbolic (resp. parabolic) elements of 
$\mbox{PSL}(2, {\mathbb R})$ are the exponentials $\exp(A)$ of \emph{hyperbolic\/} 
(resp. \emph{parabolic,\/} \emph{elliptic\/}) elements of 
${\mathcal G} = \mbox{sl}(2, {\mathbb R})$, i.e. such that 
$\mbox{det}(A) < 0$ (resp. $\mbox{det}(A)=0$, $\mbox{det}(A)>0$).

\begin{defin}
An element of $\widetilde{G}$ is \emph{hyperbolic\/} (resp. parabolic, elliptic) 
if it is the exponential of a hyperbolic (resp. parabolic, elliptic) element of 
$\mathcal G$. 
\end{defin}

\begin{defin} 
\label{def.synch}
An element $\gamma = (\gamma_{L}, \gamma_{R})$ of $\widetilde{G} \times \widetilde{G}$ 
is synchronised if, up to a permutation of left and right components, 
it has one of the following form:

\begin{itemize}

\item (hyperbolic translation): $\gamma_L$ is trivial and $\gamma_R$ is hyperbolic,

\item (parabolic translation): $\gamma_L$ is trivial and $\gamma_R$ is parabolic,

\item (hyperbolic - hyperbolic) $\gamma_L$ and $\gamma_R$ are both non-trivial and hyperbolic,

\item (parabolic - hyperbolic) $\gamma_L$ is parabolic and $\gamma_R$ is hyperbolic,

\item (parabolic - parabolic) $\gamma_L$ and $\gamma_R$ are both non-trivial and parabolic, 

\item (elliptic) $\gamma_L$ and $\gamma_R$ are elliptic elements conjugate in $\widetilde{G}$. 

\end{itemize}

An element $\gamma$ of $(\widetilde{G} \times \widetilde{G})_{/\mathcal Z}$
is synchronized if it is represented by a synchronized element of 
$\widetilde{G} \times \widetilde{G}$.

\end{defin}

\begin{lem}[Lemma $9.6$ in \cite{ba1}]
\label{lem.syncaffine}
An isometry $\gamma$ is synchronized if and only if there is an affine domain $U$ in 
$\widetilde{\mbox{AdS}}$ such that $\gamma^n(U) \cap U \neq \emptyset$ for every $n$ in 
$\mathbb Z$.
\fin
\end{lem}

Observe:

\begin{lem}[Lemma $5.6$ of \cite{ba1}]
\label{Edansds}
Every generic closed achronal subset $\Lambda$ of $\widehat{\mbox{Ein}}_{n}$ 
is contained in a 
de Sitter domain.
\fin
\end{lem}

Hence:

\begin{cor}
Any isometry preserving a generic closed achronal subset of $\widehat{\mbox{Ein}}_2$ is 
synchronized.
\fin
\end{cor}

\subsection{Causal open subsets}
\label{sub.causaldomain}

Let $\gamma = (\gamma_{L}, \gamma_{R}) = (\exp(X_L), \exp(X_R))$ be a synchronised element of 
$\widetilde{G} \times \widetilde{G}$.

\begin{defin}
The standard causal subset of $\gamma$, denoted by $C(\gamma)$, is the set of points
$x$ of $\widetilde{\mbox{AdS}}$ for which $\gamma x$ is not causally related to
$x$.
\end{defin}

Observe that $C(\gamma)=C(\gamma^{-1})$, and $C(\gamma)$ is $\gamma$-invariant.
The inclusions $C(\gamma^{k}) \subset C(\gamma)$ follows.

\begin{defin}
The convex causal subset of $\gamma$, denoted by $C_{\infty}(\gamma)$, 
is the set of points of $\widetilde{\mbox{AdS}}$ admitting in their future no 
$\gamma$-iterates of themselves.
\end{defin}

Clearly, $C_{\infty}(\gamma)$ is the decreasing intersection of all 
$C(\gamma^{n})$ when $n$ describes all $\mathbb Z$. It is $\gamma$-invariant.

At first glance, it seems natural to consider $C_{\infty}(\gamma)$ as the
prefered $\gamma$-invariant subset such  that the quotient is \emph{causal,\/} i.e. 
does not admit closed causal curves (see \cite{beem}, page $7$). Actually, 
there exists a bigger subset with the same property, which, in some way, 
is a maximal open subset with this property. The construction goes as follows: 
$\gamma$ is the time one map of the flow
$\gamma^{t}  = (\exp(tX_L), \exp(tX_R))$ induced by some Killing vector field
$X_{\gamma}$ of $\widetilde{\mbox{AdS}}$.

\begin{defin}
\label{def.absolu}
The absolute causal subset of $\gamma$, denoted by $D(\gamma)$, is the open subset 
of $\widetilde{\mbox{AdS}}$ where $X_{\gamma}$ is spacelike.
\end{defin}

\begin{lem}
The open domain $D(\gamma)$ is the union of all the $C(\gamma^{\frac{1}{n}})$.\fin
\end{lem}

Clearly, since $X_\gamma$ is lightlike on the boundary of $D(\gamma)$:

\begin{lem}
Let $U \subset \widetilde{\mbox{AdS}}$ be a $\gamma^t$-invariant subset containing $D(\gamma)$. 
If the quotient of $U$ by $\gamma$ is causal, then $U = D(\gamma)$.\fin
\end{lem}

\begin{prop}
\label{prop.Dcausal}
If $\gamma$ is not a pair $(\gamma_L, \gamma_R)$ of elliptic elements with irrationnal 
rotation angle, the quotient space of $D(\gamma)$ by $\gamma$ is a strongly
causal spacetime.
\end{prop}

\preu

Denote:

\[ R_0 = \left(\begin{array}{cc}
                 0 & 1 \\
                 -1  & 0 \end{array}\right)\]

\[ \Delta = \left(\begin{array}{cc}
                 1 & 0 \\
                 0  & -1 \end{array}\right)\]

\[ H = \left(\begin{array}{cc}
                 0 & 1 \\
                 0  & 0 \end{array}\right)\]

Up to conjugacy in $\widetilde{G} \times \widetilde{G}$, inversion of time orientation 
and permutation of the left-right components, we have 7 cases to consider:

\begin{enumerate}
\item $(X_L , X_R) = (\lambda R_0, \lambda R_0) \;\; (\lambda > 0)$, 
\item $(X_L , X_R) = (\lambda \Delta , 0 ) \;\; (\lambda > 0)$, 
\item $(X_L , X_R) = (H, 0 )$, 
\item $(X_L , X_R) = (H, -H)$, 
\item $(X_L , X_R) = (H, H)$, 
\item $(X_L , X_R) = (\lambda \Delta , \mu \Delta ) \;\; (0 < \lambda \leq \mu )$, 
\item $(X_L , X_R) = (\lambda \Delta , H ) \;\; (\lambda > 0)$.
\end{enumerate}

According to lemmas~\ref{lem.quocausal}, \ref{lem.2causal} the proposition 
is proved as soon as we check in every case that every connected 
component of $D(\gamma)$ is simply connected (when not empty). 

For every $\tilde{g}$ in $\widetilde{\mbox{AdS}} \approx \widetilde{G}$, the norm of 
$X_\gamma(\tilde{g})$ is $-det$ of $X_L - gX_Rg^{-1} = X_L - Ad(g)X_R$ (where 
$g = p(\tilde{g})$). It follows easily that in case $(3)$ and $(4)$, $D(\gamma)$ 
is actually empty.

\subsubsection{Case $(1)$: conjugacy by an elliptic element}
In this case $D(\gamma)$ is $\widetilde{G} \setminus {\bf R}$, where 
${\bf R} = \{ \exp(tR_0) \}$, i.e. the complement of the set of fixed points of 
$\gamma$. The quotient of $D(\gamma)$ by the flow $\gamma^t$ is simply
connected: apply remark~\ref{rk.quocausalfini}.

\subsubsection{Case $(2)$: translation by a hyperbolic element}
In this case the action of $\gamma$ is free and properly discontinuous since it is an action 
by left translation. $D(\gamma)$ is the entire $\widetilde{G}$: it is homeomorphic
to ${\mathbb R}^3$ hence simply connected.

\subsubsection{Case $(5)$: conjugacy by a parabolic element}
In this case $(5)$ $D(\gamma)$ is $p^{-1}(U)$, where 
$U \subset \mbox{SL}(2, {\mathbb R})$ 
is the set of matrices:

\[ \left(\begin{array}{cc}
               a & b \\
               c & d \end{array}\right) \;\;\;\; ad-bc = 1 , c \neq 0 \]

For $g$ in $U$ the iterate $\exp(nH)g\exp(-nH)$ is the matrix:

\[ \left(\begin{array}{cc}
               a+nc & -n^2c +n(d-a) +b  \\
               c & -nc + d \end{array}\right) \]

Since $c \neq 0$ it follows easily that the action on $D(\gamma)$ is free and properly
discontinuous. Every connected component of $D(\gamma)$ is simply connected.

\subsubsection{Case $(6)$: the hyperbolic-hyperbolic case}
A straightforward calculus shows that in this case $\tilde{g}$ belongs to 
$D(\gamma)$ if and only if $bc < \frac{(\lambda-\mu)^2}{4\lambda\mu}$, 
where $g= p(\tilde{g})$ is:

\[ g = \left(\begin{array}{cc}
               a & b \\
               c & d \end{array}\right) \;\;\;\; ad-bc = 1 \]

The projection in $G$ of the $\gamma^n$-iterate of $\tilde{g}$ is:

\[  \left(\begin{array}{cc}
               a\exp(n(\lambda-\mu)) & b\exp(n(\lambda+\mu)) \\
               c\exp(-n(\lambda+\mu)) & d\exp(n(\mu-\lambda)) \end{array}\right)  \]

If $\lambda \neq \mu$ then the action of $\gamma$ on the entire 
$\widetilde{\mbox{AdS}}$ is free and properly discontinuous 
(see for example \cite{salein}). The strong causality of the action
on $D(\gamma)$ is once more a corollary of lemmas~\ref{lem.quocausal}, \ref{lem.2causal}.
When $\lambda = \mu$, the projection of $D(\gamma)$ is $\{ bc < 0 \}$:
it is easy to see that the action on it is free and properly discontinuous
and we conclude thanks to lemmas~\ref{lem.quocausal}, \ref{lem.2causal}, observing
that $D(\gamma)$ is simply connected.

\subsubsection{Case $(7)$: the hyperbolic-parabolic case}
This last case is completely similar to the previous one. 
The action of $\gamma$ on 
$\widetilde{\mbox{AdS}}$ is free 
and properly discontinuous (see \cite{salein}), and $D(\gamma)$ 
is defined by:

\[ -2ac < \lambda  \]

Details are left to the reader.

\fin

\rque
In case $(6)$, if $\lambda \neq \mu$, the domain $D(\gamma)$ is not contained 
in an affine domain. This is an union of
elementary domain of invisibility $\{ bc < 0 \}$, connected by domains 
$\{ 0 \leq bc < \frac{(\lambda-\mu)^2}{4\lambda\mu} \}$. The reader can 
find complementary descriptions in \cite{brill}, or \cite{BTZ2}. 
In the terminology of these papers, $\{ bc < 0 \}$ is the union 
of regions of type I and II, and 
$\{ 0 \leq bc < \frac{(\lambda-\mu)^2}{4\lambda\mu} \}$ are the regions of type III: 
the ``inner regions''. Compare in particular our proof of 
proposition~\ref{prop.Dcausal} with {\S} $3.2.5$ of \cite{BTZ2}.
\erque

\rque 
The region $D(\gamma)$ in the case $(7)$ is particularly difficult to draw. 
The best way to catch a picture of it is to consider this case as a limit
of case $(6)$: for every $\epsilon > 0$ define $\gamma_\epsilon = (\gamma_L, \gamma_R^\epsilon)$ where:

\[ \gamma_R^\epsilon = \left(\begin{array}{cc}
             \exp(\epsilon) & \frac{\sinh(\epsilon)}{\epsilon} \\
              0 & \exp(-\epsilon)\end{array}\right)  \] 

$\gamma_\epsilon$ is the exponential of $(\lambda\Delta, \epsilon\Delta + H)$. 
At the limit $\epsilon \to 0$, $\gamma_\epsilon$ tends to $\gamma$. Then
$D(\gamma)$ is the limit of the domains 
$D(\gamma_\epsilon)$. Recall also \S~\ref{subsub.extreme}.
\erque

\rque
When $\gamma_L = \gamma_R$ are elliptic elements with rationnal angle 
the quotient of $\widetilde{\mbox{AdS}}$ by $\gamma$ is a singular spacetime 
with orbifold type.
More precisely, the timelike line of $\gamma$-fixed points induces in the quotient a singular 
line, which can be considered as the trajectory of a massive particle.

This point of view can be extended to the irrationnal angle case without difficulty, 
but we don't want to enter here in this discussion. See for example \cite{gott, matschull}.
\erque

\section{Actions on invisible domains from elementary achronal subsets}
\label{achronalinvariant}

According to Theorem $10.1$ in \cite{ba1}:

\begin{thm}
\label{noneleOK}
Let $\widetilde{\Lambda}$ be a nonelementary generic achronal subset, 
preserved by a torsionfree discrete group $\Gamma \subset \mbox{SO}_0(2,2)$. 
Then the action of
$\Gamma$ on $\Omega(\widetilde{\Lambda})$ and $E(\widetilde{\Lambda})$ 
is free, properly discontinuous and the quotient spacetime 
$M_{\widetilde{\Lambda}}(\Gamma)$ is strongly causal.
\fin
\end{thm} 

This statement does not hold when $\widetilde{\Lambda}$
is elementary. As it will appear from our study, 
the philosophy which should be retained is that in 
the elementary cases the invariant achronal subset $\widetilde{\Lambda}$, 
even if $\Gamma$-invariant is not sufficient to reveal the causal 
properties of $\Gamma$: some points are missing 
(see \S~\ref{lambdafromgamma}).

\subsection{The extreme case} 
\label{sub.invextreme}
Assume that $\widetilde{\Lambda}$ is extreme, i.e. is a lightlike segment $[x,y]$. 
We can assume that the lightlike geodesic $l$ containing $[x,y]$ 
is a leaf of the left foliation, i.e. an element of ${\mathbb R}P^1_L$. 
Let $l_x$, $l_y$ be the right leaves, i.e. the elements of ${\mathbb R}P^1_R$, 
containing respectively $x$, $y$. Then, $l$ is a fixed point of 
(the projection in $G$) of the left component $\gamma_L$ of every element of 
$\Gamma$ and $l_x$, $l_y$ are fixed points of $\gamma_R$. It follows easily 
that $\gamma$ is synchronized, that the right component $\gamma_R$ is trivial 
or hyperbolic, and that the left component is nonelliptic (maybe trivial). 
In other words, with the notations involved in the proof of 
proposition~\ref{prop.Dcausal}, $\gamma_L = \exp(\lambda\Delta + \eta{H})$ 
and $\gamma_R = \exp(\mu\Delta)$.

The first commutator group $[\Gamma, \Gamma]$ is a group of left translations.
Since $\Gamma$ is discrete, the same is true for $[\Gamma, \Gamma]$. 
Assume that $[\Gamma, \Gamma]$ is not trivial. Then it is a cyclic group 
$\mbox{Aff}(\mathbb R)$ of affine transformations of the line.
A homothety of $\mathbb R$ cannot be a commutator of elements of $\mbox{Aff}(\mathbb R)$. 
Hence, in the last case above, $[\Gamma, \Gamma] \subset \widetilde{G}_L$ is a cyclic 
group of parabolic elements preserved by conjugacies by left components of elements of 
$\Gamma$.
It follows that these left components are necessarely parabolic, i.e. translations 
of $\mathbb R$. Hence, left components of elements of $[\Gamma, \Gamma]$ 
are trivial: contradiction.

Therefore, $\Gamma$ is an abelian discrete subgroup of $A_{hyp}$, $A_{ext}$ where:

-- $A_{hyp} = \{ (\exp(\lambda\Delta), \exp(\mu\Delta)) \;\; (\lambda, \mu \in {\mathbb R})\}$,

-- $A_{ext} = \{ (\exp({\lambda}H), \exp(\mu\Delta)) \;\; (\lambda, \mu \in {\mathbb R})\}$.

\subsubsection{The mixed case $\Gamma \subset A_{ext}$}
\label{subsub.extremepar}

In this case the action is free since a parabolic element can be conjugate in 
$\widetilde{G}$ to a hyperbolic one only if they are both trivial. 

\emph{Claim: the action of $\Gamma$ on $E(\widetilde{\Lambda})$ is properly discontinuous.\/}

Let's prove now the properness: assume by contradiction the existence of a compact $K$ in 
$E(\widetilde{\Lambda})$ and a sequence $\gamma_n = (\exp({\lambda}_nH), \exp(\mu_n\Delta))$ 
of elements of $\Gamma$ such that every $\gamma_n K \cap K$ is not empty. Let 
$\Vert$ be the operator norm on $\mbox{gl}(E)$. Up to a subsequence, 
$\gamma_n/\Vert \gamma_n \Vert$ converges to an element $\bar{\gamma}$ of 
the unit ball of $\mbox{gl}(E)$. 
Since $\Gamma$ is discrete, and since all the $\gamma_n$ have determinant one, the norms 
$\Vert \gamma_n \Vert$ tends to $+\infty$. 

If the $\lambda_n$ are unbounded, up to a subsequence, we can assume that they tend to 
$+\infty$. Then, the kernel of $\bar{\gamma}$ is a hyperplane, and its image is a line.
More precisely, the image is the line spanned by one of the extremities of $[x,y]$, 
let's say $x$; and the kernel is the $Q$-orthogonal $y^\perp$. The compact $K$ is disjoint 
from $y^\perp$:
it follows that in $P(E) \setminus y^\perp$, the sequence $\gamma_n$ converges uniformly on 
$K$ towards the constant map $x$. This is a contradiction, since $x$ does not belong to $K$.

If the $\lambda_n$ are bounded, the image and the kernel of $\bar{\gamma}$ are $Q$-isotropic 
$2$-planes (one of them is the $2$-plane spanned by $x$ and $y$): their projection in 
$S(E)$ is disjoint from ${\mathbb A}{\mathbb D}{\mathbb S}$.
But the iterates $\gamma_n K$ accumulates on the projection of the image of $\bar{\gamma}$:
we obtain a contradiction as above. The claim is proved.

Moreover, according to proposition~\ref{prop.Dcausal}, case $(7)$, this action is 
strongly causal 
except if some right component $\exp(\mu\Delta)$ is trivial (case $(3)$ of 
proposition~\ref{prop.Dcausal}): in this last case, the quotient is 
foliated by closed lightlike geodesics, which are orbits of some 
$1$-parameter subgroup of $A_{ext}$.

\subsubsection{The hyperbolic case $\Gamma \subset A_{hyp}$} 
\label{subsub.extremehyp}
There is a particular situation: the subcase $\Gamma \subset \widetilde{G}_L$. 
Then, $\Gamma$ is cyclic. It follows from proposition~\ref{prop.Dcausal}, case $(2)$, 
that the action on $E(\widetilde{\Lambda})$ is free, properly discontinuous and 
strongly causal. 
The same conclusion holds if $\Gamma \subset \widetilde{G}_R$.

Hence, assume that $\Gamma$ is not contained in $\widetilde{G}_R$ or $\widetilde{G}_L$.
The group $A_{hyp}$ admits $4$ fixed points in $\overline{\mbox{Ein}}_2$, including the 
projections of $x$, $y$. We can then define two additionnal $A_{hyp}$-fixed points 
$x'$, $y'$ uniquely defined by the requirement that $\{ x,x' \}$ and 
$\{ y,y' \}$ are strictly achronal. 

Many subcases appear, with different behavior. For example, the action of $\Gamma$ 
on $E(\widetilde{\Lambda})$ may be free and properly discontinuous (for example, 
if $\Gamma$ is cyclic, spanned by an element for which $\lambda > \mu$). 
But the action may also be non proper (the cyclic case, with $\lambda = \mu$). 
Anyway, this action is never causal. Indeed, if $\gamma$ is an element of 
$\Gamma \setminus (\widetilde{G}_L \cup \widetilde{G}_R)$, $E(\widetilde{\Lambda})$ 
is $\gamma^t$-invariant, but is not contained in the absolute causal domain 
$D(\gamma)$. Then, the $\gamma^t$-orbit of a point $x$ in 
$E(\widetilde{\Lambda}) \setminus D(\gamma)$ is a timelike curve 
containing $x$ and $\gamma x$.

\subsection{The splitting case}
\label{splitcausal}
We consider the splitting case $\widetilde{\Lambda} = \{ x,y \}$, with $x$, $y$ 
not causally related. Then, the leaves of $\widehat{\mathcal G}_R$ through $x$, 
$y$ are two distinct fixed points in ${\mathbb R}P^1_R$. The right component of 
any element of $\Gamma$ is therefore trivial or hyperbolic. A similar argument 
shows that the left components are trivial or hyperbolic too. Hence, 
in the notations of the previous {\S}, we have $\Gamma \subset A_{hyp}$.

Observe that the segment $]x,y[$ is contained in $E(\widetilde{\Lambda})$. 
Hence, if $\Gamma$ is not cyclic, its action on $E(\widetilde{\Lambda})$ 
cannot be properly discontinuous. 

Assume that $\Gamma$ is cyclic, spanned by $\gamma = (\exp(\lambda\Delta), \exp(\mu\Delta))$. 
If $\lambda$ or $\mu$ is zero the action is free, properly discontinuous and causal. 

If $\lambda$ and $\mu$ are both nonzero, it follows from case $(6)$ of 
proposition~\ref{prop.Dcausal} that \emph{the action of $\Gamma$ on 
$E(\widetilde{\Lambda})$ is free, properly discontinuous and strongly causal 
if and only if $x$, $y$ are attractive or repulsive fixed points of $\gamma$.\/}

\subsection{The conical case} 
\label{conicalbof}
We assume here that $\widetilde{\Lambda}$ is conical, i.e. 
the union of two non-trivial lightlike segments $I_1 = [x_1, x]$ and 
$I_2 = [x, x_2]$. 
Then, $\{ x_1, x_2 \}$ is strictly achronal, $E(\widetilde{\Lambda})$ 
is contained in $E(x_1, x_2)$, and $\Gamma$ preserves $E(x_1, x_2)$. 
As in the previous {\S} we still have $\Gamma \subset A_{hyp}$.

Recall the description of $E(\widetilde{\Lambda})$ (\S~\ref{subsub.conical}): 
it is the union of
a outer region $P_2 \cap F_2$, an intermediate region $F_2 \cap F_1$, 
and their common horizon boundary. It follows that $E(\widetilde{\Lambda})$ 
is $A_{hyp}$-invariant, and that all the $A_{hyp}$-orbits inside 
$E(\widetilde{\Lambda})$ are $2$-dimensional. Hence, the action of 
$A_{hyp}$ on $E(\widetilde{\Lambda})$ is free and properly discontinuous: 
the same is true 
for the action of $\Gamma$.

If $\Gamma$ is contained in $\widetilde{G}_L$ or $\widetilde{G}_R$, then its action on 
$E(\widetilde{\Lambda})$ is strongly causal. If not, the statement in the previous case still holds: 
the action of $\Gamma$ on $E(\widetilde{\Lambda})$ is free, properly discontinuous and 
strongly causal  
if and only if $\Gamma$ is cyclic and $x$, $y$ are attractive or repulsive fixed points of 
every non-trivial element of $\Gamma$.

\section{Existence of invariant achronal subsets}
\label{sec.abadmissible}

Recall the following definition (definition $10.6$ in \cite{ba1}):

\begin{defin}
\label{def.admissible}
Let $\rho_L: \Gamma \rightarrow G$ and $\rho_R: \Gamma \rightarrow G$ 
two morphisms.
The representation $\rho = (\rho_L, \rho_R)$ is admissible if and only 
if it is faithfull, has discrete image and lifts to  some representation 
$\tilde{\rho}: \Gamma \rightarrow (\widetilde{G} \times 
\widetilde{G})_{/{\mathcal Z}}$ preserving a generic closed achronal subset 
of $\widehat{\mbox{Ein}}_2$ containing at least two points.

A $\rho$-admissible closed subset for an admissible representation 
$\rho$ is the projection in $\overline{\mbox{Ein}}_2$ of 
$\tilde{\rho}$-invariant generic closed achronal subset of 
$\widehat{\mbox{Ein}}_2$ containing at least two points.
\end{defin}

In \cite{ba1}, we claimed (theorem $10.7$):

\begin{thm}
\label{thm.admissible}
Let $\Gamma$ be a torsionfree group, and 
$\rho: \Gamma \rightarrow G \times G$ a faithfull representation.
Then, $\rho$ is admissible if and only if one the following occurs:
\begin{enumerate}

\item \emph{The abelian case:\/} $\rho(\Gamma)$ is a discrete subgroup of $A_{hyp}$, 
$A_{ext}$ or $A_{par}$ where (see the notations in {\S}~\ref{sub.invextreme} 
where the first two groups are already defined):

-- $A_{hyp} =\{ (\exp(\lambda\Delta), \exp(\mu\Delta)) / \lambda, \mu \in {\mathbb R}\}$,

-- $A_{ext} =\{ (\exp(\lambda\Delta), \exp(\eta{H})) / \lambda, \eta \in {\mathbb R}\}$,

-- $A_{par} =\{ (\exp(\lambda{H}), \exp(\lambda{H})) / \lambda \in {\mathbb R}\}$.

\item \emph{The non-abelian case:\/} The left and right morphisms 
$\rho_L$, $\rho_R$ are faithfull with discrete image and the marked surfaces 
$\Sigma_L = \rho_L(\Gamma)\backslash{\mathbb H}^2$, 
$\Sigma_R = \rho_R(\Gamma)\backslash{\mathbb H}^2$ are homeomorphic, 
i.e. there is a $\Gamma$-equivariant homeomorphism 
$f: {\mathbb H}^2 \rightarrow {\mathbb H}^2$ satisfying:

\[ \forall \gamma \in \Gamma, \; f \circ \rho_L(\gamma) = \rho_R(\gamma) \circ f\]

\end{enumerate}
\end{thm}

\rque
\label{rk.endpara}
Parabolic $\rho_L(\gamma)$ correspond to punctures in $\Sigma_L$. Hence, 
if $\gamma$ is a homotopy class corresponding to a loop which is not homotopic 
to an isolated end of $\Sigma_R$, $\rho_L(\gamma)$ is necessarely hyperbolic. 
\erque

Actually, we only proved in \cite{ba1} the non-abelian case. Here, we
justify the abelian case.

When $\Gamma$ is cyclic, a case by case study is needed, 
but which follows almost immediatly from the study in the 
proof of Proposition~\ref{prop.Dcausal}. 
The situation can be summarized as follows:

\begin{prop}
\label{pro.cyclicadmissible}
When $\Gamma$ is cyclic, then the representation $\rho$ 
is admissible if and only if either the left or right component 
of $\rho(\gamma)$ is hyperbolic and the other component non-elliptic, 
or if $\rho_L(\gamma)$, $\rho_R(\gamma)$ are parabolic elements conjugate 
one to the other in $G$.\fin
\end{prop}

Assume now that $\Gamma$ is abelian but not cyclic. 
It follows from the cyclic case that $\rho$ is 
admissible if and only if $\rho(\Gamma)$ is
contained in (a conjugate of $G$ of) $A_{hyp}$, $A_{ext}$ or $A_{par}$. 
The validity of Theorem~\ref{thm.admissible} in the abelian case follows.

\section{Minimal invariant achronal subsets}
\label{lambdafromgamma}
Let $\rho: \Gamma \rightarrow G \times G$ be an admissible representation. 

\begin{defin}
$\overline{\Lambda}(\rho)$ is the closure of the set of attractive fixed points 
in $P(E)$.
\end{defin}

Since attractive fixed points in $P(E)$ of elements of $G$ belong to 
$\overline{\mbox{Ein}}_2$, $\overline{\Lambda}(\rho)$ is contained in 
$\overline{\mbox{Ein}}_2$.
According to {\S} $10.5$ of \cite{ba1}:

\begin{thm}
\label{thm.minimal}
Let $\Gamma$ be a non-abelian torsion-free group, and 
$\rho: \Gamma \rightarrow G \times G$.
Then, every $\rho(\Gamma)$-invariant closed subset of 
$P(E)$ contains $\overline{\Lambda}(\rho)$.
\end{thm}

\begin{cor}
\label{cor.cor}
Let $(\Gamma,\rho)$ be pair satisfying the hypothesis of Theorem~\ref{thm.minimal}. Then,
$\overline{\Lambda}(\rho)$ is a $\rho(\Gamma)$-invariant generic nonelementary achronal 
subset of 
$\overline{\mbox{Ein}}_2$. Furthermore, for every $\rho(\Gamma)$-invariant closed achronal 
subset 
$\Lambda$ in $\mbox{Ein}_2$, the invisibility domain $E(\Lambda)$ projects injectively in 
$\overline{{\mathbb A}{\mathbb D}{\mathbb S}}$ inside $E(\overline{\Lambda}(\rho))$.
\fin
\end{cor}

\rque
\label{le.dansS}
We will need the following remark: still assuming that $\Gamma$ is not abelian,
the minimal closed achronal subset $\overline{\Lambda}(\rho)$ is the projection
in $P(E)$ of $\Lambda(\rho) \cup -\Lambda(\rho)$, where $\Lambda(\rho)$ is a closed
achronal subset of $\mbox{Ein}_2$, and $-\Lambda(\rho)$ is the image of $\Lambda(\rho)$
by the antipody in $S(E)$. The action of $\rho(\Gamma)$ on $\Lambda(\rho)$ is minimal,
hence $\Lambda(\rho)$ and $-\Lambda(\rho)$ are exactly the minimal components of
the action of $\rho(\Gamma)$ on $S(E)$ (see Lemma $10.21$ in \cite{ba1}).
\erque

When $\Gamma$ is abelian, Theorem~\ref{thm.minimal} and Corollary~\ref{cor.cor} 
do not hold. However, in this case, still assuming that $\rho: \Gamma \rightarrow G$ 
is admissible:

-- \emph{The extreme case:\/} let $l$ be the unique $A_{ext}$-invariant left leaf, and $r_1$, 
$r_2$ be the two $A_{ext}$ right invariant right leaves.  Let $\overline{\Lambda}_{ext}$ 
be the set of $A_{ext}$-fixed points: $\overline{\Lambda}_{ext} = \{ l \cap r_1, l\cap r_2\}$. 
It is easy to show that if $\rho(\Gamma) \subset A_{ext}$, any generic $\rho(\Gamma)$-invariant 
closed achronal subset containing at least two points  must contain $\overline{\Lambda}_{ext}$.
This formulation is an extension - more accurately, a limit case - of corollary~\ref{cor.cor}, 
even if elements of $\rho(\Gamma)$ do not admit attractive fixed points in $P(E)$.

-- \emph{The parabolic case:\/} if $\rho(\Gamma) \subset A_{par}$, any $\rho(\Gamma)$-invariant 
closed subset contains the unique fixed point of $A_{par}$. Hence, for any 
$\rho(\Gamma)$-invariant closed achronal subset $\Lambda$, we have: 
$E(\Lambda) \subset D(\rho(\Gamma))$. This is in some way a limit case of the previous one.

-- \emph{The hyperbolic-hyperbolic case:\/} the group $A_{hyp}$ admits $4$ fixed points 
in $P(E)$. If $\rho(\Gamma)$ is a lattice of $A_{hyp}$, it is easy to see that every 
$A_{hyp}$-fixed point is an attractive fixed point in $P(E)$ of some $\rho(\gamma)$. 
Anyway, corollary~\ref{cor.cor} is false in this situation. There is however a 
convenient statement:
any nonelementary $\rho(\Gamma)$-invariant generic closed achronal subset 
of $P(E)$ must 
contain the $4$ fixed points of $A_{hyp}$.

\section{Coincidence of standard and absolute chronological domains}
\label{sec.coincide}
In \S~\ref{sub.causaldomain}, we have associated to any synchronized element $\gamma$ of $\widetilde{G} \times \widetilde{G}$ two open domains:

-- the convex causal domain $C_{\infty}(\gamma)$,

-- the absolute causal subset $D(\gamma)$.

In all cases, we have $C_{\infty}(\gamma) \subset D(\gamma)$. Actually, the identity $C_\infty(\gamma) = D(\gamma)$ holds if and only if $\gamma$ is a pair $(\gamma_L,\gamma_R)$ where $\gamma_L$, $\gamma_R$ are conjugate in $\widetilde{G}$.

These definitions easily extend to (lifted) admissible representations $\widetilde{\rho}: \Gamma \rightarrow 
\widetilde{G} \times \widetilde{G}$:

\begin{defin}
\label{def.causaldomains}
The convex causal domain $C_\infty(\widetilde{\rho})$ is the interior of the intersection $\bigcap_{\gamma \in \Gamma} C_\infty(\widetilde{\rho}(\gamma))$.
The absolute causal domain $D(\widetilde{\rho})$ is the interior of the intersection $\bigcap_{\gamma \in \Gamma} D(\widetilde{\rho}(\gamma))$.
\end{defin}

The inclusions $C_{\infty}(\gamma) \subset D(\gamma)$ implies $C_{\infty}(\widetilde{\rho}) \subset D(\widetilde{\rho})$. Conversely:

\begin{thm}
\label{thm.egal}
If $\Gamma$ is non cyclic, then the convex causal domains and absolute causal domains coincide.
\end{thm}

The rest of this section is devoted to the proof of theorem~\ref{thm.egal}.

\subsection{The flat case}
\label{sub.thmflat}
In the flat case, i.e., the case where $\widetilde{\rho}(\Gamma)$ preserves a point in $\widetilde{\mbox{AdS}}$, the proof of Theorem~\ref{thm.egal} is obvious. Indeed, after conjugacy, we can assume in this case that the left and right representations $\rho_L$, $\rho_R$ coincide. Then, for every $\gamma$ in $\Gamma$, the identity $C_{\infty}(\rho(\gamma)) = D(\gamma)$ holds.

\subsection{The abelian case}
\label{sub.thmabelian}
If $\Gamma$ is abelian, since it is assumed non-cyclic, $\rho(\Gamma)$ is a lattice in $A_{hyp}$ or $A_{ext}$. We only consider the first case, the other can be obtained in a similar way (or as a limit case). 

There are two morphisms $\alpha,\beta: \Gamma \rightarrow \mathbb R$ such that, for every $\gamma$ in $\Gamma$:
\[ \rho_L(\gamma) = \left(\begin{array}{cc}
                                                                                  \exp(\alpha(\gamma)) & 0 \\
                                                                                  0 & \exp(-\alpha(\gamma)) 
                                                                                  \end{array}\right), \;\;
\rho_R(\gamma) =                                               \left(\begin{array}{cc}
                                                                                  \exp(\beta(\gamma)) & 0 \\
                                                                                  0 & \exp(-\beta(\gamma)) 
                                                                                  \end{array}\right) \]

Recall that the projections in ${\mathbb A}{\mathbb D}{\mathbb S} \approx \mbox{SL}(2, {\mathbb R})$ are:

\[ C_\infty(\rho(\gamma)) = \{ bc < 0 \}, \;\;\; D(\rho(\gamma)) = \{ bc < \sinh^2(\alpha(\gamma)- \beta(\gamma)) \} \]

Theorem~\ref{thm.egal} follows then from the fact that, since $\rho(\Gamma)$ is a lattice in $A_{hyp}$,
$\mid \alpha(\gamma) - \beta(\gamma) \mid$ admits arbitrarly small value.

\subsection{The proper case}
Assume that $\rho$ is strongly irreducible. The representations $\rho_L$ and $\rho_R$ are faithfull, 
with discrete image, semi-conjugate one to the other, but not conjugate in $G$, and $\Gamma$ is not abelian. 

Let $\Lambda(\rho)$ be one of the two minimal closed subsets of ${\mathcal D} \subset S(E)$ such 
that the closure in $S(E)$ of the set of attractive fixed points is $\Lambda(\rho) \cup -\Lambda(\rho)$ 
(see remark~\ref{le.dansS}). It projects injectively on $\overline{\Lambda}(\rho)$, the closure in $P(E)$ 
of the set of attractive fixed points of elements of $\rho(\Gamma)$.

Let $E(\rho)$ be the invisible domain of $\Lambda(\rho)$: this is the intersection between ${\mathbb A}{\mathbb D}{\mathbb S}$ and the intersection of all $\{ x / \langle x \mid p \rangle < 0 \}$, where $p$ describe $\Lambda(\rho)$. Let $\overline{E}(\rho)$ be the projection of $E(\rho)$ in $\overline{{\mathbb A}{\mathbb D}{\mathbb S}}$. It can be defined in the following way:

\[ \overline{E}(\rho) = \{ [x] \in \overline{{\mathbb A}{\mathbb D}{\mathbb S}} / 
\forall p, q \in {\Lambda}(\rho), \langle x \mid p \rangle\langle x \mid q \rangle > 0\}\]

Indeed, although $\langle x \mid p \rangle$, $\langle x \mid q \rangle$ are not individually well 
defined for $[x]$ in $P(E)$, their product have a well-defined sign.

For any subset $J$ of $\overline{\Lambda}(\rho) \approx {\Lambda}(\rho)$, we can define 
$\overline{E}(J)$ as the interior of the set: 
$\{ [x] \in \overline{{\mathbb A}{\mathbb D}{\mathbb S}} / \forall [p], [q] \in J, 
\langle x \mid p \rangle\langle x \mid q \rangle > 0\}$.

\begin{lem}
\label{le.dense}
If $J$ is $\rho(\Gamma)$-invariant and non-empty, then $\overline{E}(J) = \overline{E}(\rho)$.
\end{lem}

\preu
The inclusion $\overline{E}(\rho) \subset \overline{E}(J)$ is obvious (observing that $\overline{E}(\rho)$ is open). The reverse inclusion follows from the fact that $\rho(\Gamma)$-invariant subsets of $\overline{\Lambda}(\rho)$ are dense, that $\overline{E}(J)$ is open if $J$ is closed, and that if $\bar{J}$ is the closure of $J$, $\overline{E}(\bar{J}) = \overline{E}(J)$.\fin

\begin{cor}
\label{cor.identifionsC}
Every connected component of $C_\infty(\widetilde{\rho})$ projects injectively in $P(E)$ on 
$\overline{E}(\rho)$.
\end{cor}

\preu
Let $\gamma$ be a non-trivial element of $\Gamma$. We define $J(\gamma) \subset \overline{\mathcal Q}$ in the following way:

-- if $\rho(\gamma)$ is hyperbolic - hyperbolic: $J(\gamma) =\{ p^+(\gamma),p^-(\gamma) \}$, where $p^+(\gamma)$ is the attractive fixed point of $\rho(\gamma)$, and $p^-(\gamma)$ is the repulsive fixed point.

-- if $\rho(\gamma)$ is hyperbolic - parabolic: $J(\gamma) = \{ p(\gamma), q(\gamma) \}$, where $p(\gamma)$, $q(\gamma)$ are the two $\rho(\gamma)$ fixed points,

-- if $\rho(\gamma)$ is parabolic - parabolic: $J(\gamma) = \{ p(\gamma) \}$,
where $p(\gamma)$ is the unique fixed point.

Then, the study in the proof of Proposition~\ref{prop.Dcausal} shows that, in every case,
every connected component of $C_\infty(\rho(\gamma))$ projects injectively on $\overline{C}_\infty(\rho(\gamma)) = \overline{E}(J(\gamma))$. On the other hand, if $\rho_L(\gamma)$ (or $\rho_R(\gamma)$) is parabolic, every closed subset of ${\mathbb R}P^1_L$ (or ${\mathbb R}P^1_R$) which is $\rho_L(\gamma)$-invariant (or $\rho_R(\gamma)$-invariant) contains the unique $\rho_L(\gamma)$-fixed point (resp. the $\rho_R(\gamma)$-fixed point). Hence, we have the inclusion $\overline{J}(\gamma) \subset \overline{\Lambda}(\rho)$. In other words, $\overline{C}_\infty(\rho(\gamma))$ is contained in $\overline{E}(\rho)$. The corollary~\ref{cor.identifionsC} follows then from lemma~\ref{le.dense}.\fin

\begin{cor}
\label{cor.unsuffit}
Let $\gamma_1$ be an element of $\Gamma$ such that $\rho_L(\gamma_1)$ and $\rho_R(\gamma_1)$ are both hyperbolic. Then: $C_{\infty}(\widetilde{\rho}) = \bigcap_{\gamma \in \Gamma} C_{\infty}(\widetilde{\rho}(\gamma\gamma_1\gamma^{-1}))$.
\end{cor}

\preu
Corollary of lemma~\ref{le.dense} and corollary~\ref{cor.identifionsC}, since $\bigcap_{\gamma \in \Gamma} \overline{C}_{\infty}({\rho}(\gamma\gamma_1\gamma^{-1}))$ is equal to $\overline{E}(J)$, where $J$ is the $\Gamma$-orbit of the attractive fixed points of $\rho(\gamma_1)$ and $\rho(\gamma_1^{-1})$.\fin

\preud{thm.egal}
According to \S~\ref{sub.thmflat} and \ref{sub.thmabelian}, we just have to consider the case where $\rho(\Gamma)$ is non-abelian and does not preserve a point in AdS.
The surfaces $\Sigma_L = \rho_L(\Gamma)\backslash{\mathbb H}^2$ and $\Sigma_R = \rho_R(\Gamma)\backslash{\mathbb H}^2$ are homeomorphic: let $\Sigma$ be any surface homeomorphic to $\Sigma_R$, $\Sigma_L$.

Let $c_1$ be a closed loop in $\Sigma$ which is not freely homotopic to an isolated end of $\Sigma$. It represents a conjugacy class $[\gamma_1]$ in $\Gamma$. 
According to remark~\ref{rk.endpara}, every $\rho_{L,R}(\gamma_1)$ is hyperbolic.
After conjugacy, we can assume:
\[ \rho_L(\gamma_1) = \left(\begin{array}{cc}
               \exp(\lambda) & 0 \\
               0 & \exp(-\lambda) \end{array}\right), \;\;\; \rho_R(\gamma_1) = \left(\begin{array}{cc}
               \exp(\mu) & 0 \\
               0 & \exp(-\mu) \end{array}\right) \]
with $\lambda \geq \mu>0$.

Since $\Gamma$ is non-abelian, the Euler characteristic of $\Sigma$ is negative. Hence, there is a closed loop $c_2$ in $\Sigma$,  not freely homotopic to an end of $\Sigma$, and such that every loop freely homotopic to $c_1$ intersects every loop freely homotopic to $c_2$.

Let $\gamma_2$ be any element of $\Gamma$ corresponding to the free homotopy class of $c_2$. 
We express the coefficients of $\rho_{L,R}(\gamma_2) = \exp(A_{L,R})$:
\[ A_L = \left(\begin{array}{cc}
               \alpha_L & \beta_L \\
               \nu_L & -\alpha_L \end{array}\right), \;\;\; A_R = \left(\begin{array}{cc}
               \alpha_R & \beta_R \\
               \nu_R & -\alpha_R \end{array}\right) \]

The fixed points in ${\mathbb R}P^1_L$ of $\rho_L(\gamma_1)$ are $0$ and $\infty$.
Hence, the connected components of the complement in ${\mathbb R}P^1_L$ of these fixed points are $]-\infty,0[$ and $]0,+\infty[$.

The fixed points in ${\mathbb R}P^1_L$ of $A_L$ are $\alpha_L\frac{1 \pm \sqrt{1+\beta_L\nu_L/\alpha_L^2}}{\nu_L}$. Replacing $\gamma_2$ by its inverse if necessary, and since the intersection number between $c_1$ and $c_2$ is not trivial,
we can assume that the attractive $\rho_L(\gamma_2)$-fixed point belongs to $]0,+\infty[$,
and that the repulsive fixed point belongs to $]-\infty,0[$. In other words, we can assume
that the products $\beta_L\nu_L$ and $\alpha_L\nu_L$ are positive. Hence, after conjugacy by a diagonal matrix, we can assume $\nu_{L}=\beta_{L} \neq 0$.

Then, since $\rho_L$ and $\rho_R$ are semi-conjugate, the right components satisfy the same properties: we can assume $\beta_R=\nu_R \neq 0$.

Assume now by contradiction that the inclusion $C_{\infty}(\widetilde{\rho}) \subset D(\widetilde{\rho})$ is strict. Then, there is an element $\tilde{x}$ of $D(\widetilde{\rho})$ which is in the boundary of $C_{\infty}(\widetilde{\rho})$. Let $\widetilde{U}$ be an open neighborhood of $\tilde{x}$ in $D(\widetilde{\rho})$. Let $U$ be the projection of $\widetilde{U}$ in $P(E)$. Then, according to corollary~\ref{cor.unsuffit}, there is an element $\gamma$ of $\Gamma$, a fixed point $x_1$ of $\rho(\gamma_1)$, and an element $p$ of $U$ such that 
$\langle \rho(\gamma)x_1 \mid p \rangle = 0$. After conjugacy of $\gamma_1$ by $\gamma$, we can assume that $\gamma$ is trivial. Moreover, we can also assume without loss of generality that $x_1$ is the attractive fixed point of $\rho(\gamma_1)$. 
Then, the equation $\langle x_1 \mid p \rangle = 0$ means $c=0$ where $p$ is expressed by the matrix:

\[ g =\left(\begin{array}{cc}
         a & b \\
         c & d \end{array}\right) \]

Observe that we can also assume, after a slight modification of $p \approx g$ if necessary, $b \neq 0$.

We consider the elements $\gamma_n = \gamma_1^n\gamma_2\gamma^{-n}_1$ of $\Gamma$.
We have: 

\begin{eqnarray*}
\rho_L(\gamma_n)  =  \exp(n\lambda\Delta)\exp(A_L)\exp(-n\lambda\Delta) \\
\rho_R(\gamma_n)  =  \exp(n\mu\Delta)\exp(A_R)\exp(-n\mu\Delta)
\end{eqnarray*}

Hence, the norm at $g$ of the Killing vector field generating $\gamma_n$ is the opposite of the determinant of $X_n$, with:

\[ X_n =  \exp(n\lambda\Delta)\exp(A_L)\exp(-n\lambda\Delta)g - g\exp(n\mu\Delta)\exp(A_R)\exp(-n\mu\Delta) \]

After computation, we see that $X_n$ is the matrix:
\begin{small}
\[ \left(\begin{array}{cc}
 a(\alpha_L - \alpha_R) - b\beta_R\exp(-2n\mu) & b(\alpha_L+\alpha_R) + d\beta_L\exp(2n\lambda)-a\beta_R\exp(2n\mu) \\
a\beta_L\exp(-2n\lambda) - d\beta_R\exp(-2n\mu) & b\beta_L\exp(-2n\lambda) - d(\alpha_L-\alpha_R) \end{array}\right)\]
\end{small}

We distinguish two subcases:

\subsubsection{The case $\lambda > \mu$:}
Observe that $b$ and $d$ are nonzero. If $\lambda > \mu$, the leading term of $-\mbox{det}(X_n)$ for $n \to+\infty$ is:
$$-d^2\beta_L\beta_R\exp(2n(\lambda-\mu))$$ 
On the other hand, the leading term for $n \to -\infty$ is $b^2\beta_L\beta_R\exp(-2n(\lambda+\mu))$.
But, since $g$ correspond to an element of $D(\widetilde{\rho})$, all the $-\mbox{det}(X_n)$ are positive. Hence, the product $\beta_L\beta_R$ must be positive and negative. Contradiction.\fin

\subsubsection{The case $\lambda = \mu$:}
More precisely, the remaining case is the case where $\lambda=\mu$ for any choice of pairs $\gamma_1$, $\gamma_2$ as above, such that the corresponding homotopy classes have non-trivial intersection number. It is equivalent to the fact that $\mbox{Tr}(\rho_L(\gamma_1)) = \mbox{Tr}(\rho_R(\gamma_1))$ for every $\gamma_1$ in $\Gamma$ which is not freely homotopic to a loop around an isolated end of $\Sigma$.

Select such a pair $(\gamma_1, \gamma_2)$ of $\Gamma$ satisfying the following additionnal property: the product $\gamma_3 =\gamma_1\gamma_2$ is not freely homotopic to an isolated end of $\Sigma$ (we leave to the reader the proof of the fact that such a pair exists). Let $\Gamma_1$ be the group generated by $\gamma_1$, $\gamma_2$. Hence, we can assume the identity $\mbox{Tr}(\rho_L(\gamma_i)) = \mbox{Tr}(\rho_R(\gamma_i))$ for $i=1,2,3$. By Fricke-Klein Theorem (\cite{fricke, klein, goldman}) these equalities imply that the restrictions of $\rho_L$and $\rho_R$ to $\Gamma_1$ are representations conjugated in $\mbox{SL}(2,{\mathbb R})$.

Therefore, after conjugacy,we can assume: $\rho_L(\gamma_i) = \rho_R(\gamma_i) \;\; (i=1,2)$. 
In other words, $\lambda = \mu$, $\alpha_L = \alpha_R$, $\beta_L = \beta_R$. A straightforward computation shows that the leading term of $-\mbox{det}(X_n)$ for $n \to +\infty$ is $(a\beta_L - d\beta_R)(d\beta_L - a\beta_R) = -(a-d)^2\beta_L^2$. We obtain a contradiction since this term should be nonnegative, whereas $d = 1/a \neq a$.\fin

\section{Conformal boundaries of strongly causal spacetimes}
\label{sub.bord}
As we have seen in the introduction, the notion of black-hole is related to 
the notion of conformal boundary. 

\begin{defin}[Compare with {\S} $4.2$ in \cite{bordads}]
\label{def.bord}
An {AdS-spacetime with boundary\/} is a triple $(M, O, {\mathcal M})$ where 
$\mathcal M$ is a manifold with boundary $O$ and interior $M$ which is 
$(\mbox{AdS}, \mbox{Ein}_2)$-modeled, i.e.:

-- there exist a morphism (the holonomy representation) 
$\rho = (\rho_L,\rho_R): \Gamma \rightarrow G \times G$, where 
$\Gamma$ is the fundamental group of $\mathcal M$,

-- there exist a $\rho$-equivariant local homeomorphism (the developing map) 
${\mathcal D}: \widetilde{\mathcal M} \rightarrow \overline{\mbox{Ein}}_3$, 
where $\widetilde{\mathcal M}$ is the universal covering of $\overline{M}$,

-- the image by ${\mathcal D}$ of $\widetilde{M}$, the interior of 
$\widetilde{\mathcal M}$, is contained in 
$\overline{{\mathbb A}{\mathbb D}{\mathbb S}}$,

-- the image by ${\mathcal D}$ of the boundary $\widetilde{O}$ of 
$\widetilde{\mathcal M}$ is contained in $\overline{\mbox{Ein}}_2$.

\end{defin}

Observe that if $(M, O, {\mathcal M})$ is AdS-spacetime with boundary, 
$M$ inherits a well-defined AdS-structure, and $O$ a $\mbox{Ein}_2$-structure.

Our aim is to attach to any AdS-spacetime $M$ a AdS-spacetime with boundary 
$(M, O, {\mathcal M})$. This procedure should be canonical.

\begin{defin}
Let $(M, O, {\mathcal M})$, $(M', O', {\mathcal M}')$ be 
two AdS-spacetime with 
boundary. 
A morphism between them 
is a local  homeomorphism from $\mathcal M$ into ${\mathcal M}'$ 
inducing a AdS-morphism $M \to M'$.
\end{defin}

Such a morphism lifts to a map $\tilde{F}: \widetilde{\mathcal M} \rightarrow 
\widetilde{\mathcal M}'$ such that ${\mathcal D}' \circ \tilde{F} = g \circ 
{\mathcal D}$ for some isometry $g$ of $\overline{{\mathbb A}{\mathbb D}{\mathbb S}}$. 
In particular, it induces a $\mbox{Ein}_2$-morphism $O \to O'$.
Such a morphism if an \emph{isomorphism\/} if it is moreover a homeomorphism.

\begin{defin}
An AdS-spacetime with boundary $(M, O, {\mathcal M})$ is an universal conformal 
completion of $M$ if for any AdS-spacetime with boundary 
$(M, O', {\mathcal M}')$
there exist an injective morphism $(M, O', {\mathcal M}') \to (M, O, {\mathcal M})$.
\end{defin}

It should be clear to the reader that if a AdS-spacetime $M$ admits an universal 
conformal completion, then this completion is unique up to isomorphism. 
In this case, the boundary is denoted by $O_M$, and called the \emph{natural conformal 
boundary of $M$.}

C. Frances proved that \emph{complete\/} AdS spacetimes, i.e. quotients 
of the entire AdS by  discrete torsion-free subgroups, 
admits an universal completion (\cite{bordads}, Theorem $1$).
But our spacetimes are \emph{never} complete, and we will see that
some of them do not admit universal conformal completion as defined
above (see remark~\ref{rk.nocompletion}). However, we can prove the existence
of such universal completions if we restrict to the strongly causal category:

Observe first that causal curves in an AdS-spacetime with boundary is a well defined 
notion since they are well defined in $\overline{{\mathbb A}{\mathbb D}{\mathbb S}} 
\cup \overline{\mbox{Ein}}_2$. We can therefore define the causality
relation in such a manifold with boundary, and in particular the strong
causality property (see \S~\ref{sec.strong}).

\begin{defin}
An AdS-spacetime with boundary $(M, O, {\mathcal M})$ is an universal 
strongly causal conformal 
completion of $M$ if it is strongly causal and 
for any strongly causal AdS-spacetime with boundary 
$(M, O', {\mathcal M}')$
there exist an injective morphism $(M, O', {\mathcal M}') \to (M, O, {\mathcal M})$.
\end{defin}

Obviously, an AdS-spacetime can admits a conformal strongly causal completion
only if it is already strongly causal. Conversely:

\begin{thm}
\label{thm.bordcausal}
Every strongly causal AdS-spacetime admits an universal strongly causal 
conformal completion.
\end{thm}

In the proof we will need the following lemmas, valid for any local homeomorphism 
$\varphi: X \rightarrow Y$ between manifolds (for proofs, see for example \cite{barHLP}, 
{\S} $2.1$):

\begin{lem}[Lemme des assiettes embo\^{\i}t\'ees]
\label{lem.fermassiette}
Let $U$, $U'$ be two open domains in $X$ such that the restrictions 
of $\varphi$ on $U$, $U'$ are injective. Assume that $U \cap U'$ 
is not empty, and that $\varphi(U')$ contains $\varphi(U)$. Then, 
$U'$ contains $U$.
\fin
\end{lem}

\begin{lem}[Fermeture des assiettes]
\label{lem.ferme}
Assume that $\varphi$ is injective on some open domain $U$ in $X$, 
and that the image $V = \varphi(U)$ is locally connected in $Y$, i.e. 
that every point $y$ in the closure of $V$ admits arbitrarly small 
neighborhood $W$ such that $V \cap W$ is connected. Then, the restriction 
of $\varphi$ to the closure of $U$ in $X$ is injective.
\fin
\end{lem}

\preud{thm.bordcausal}

\emph{Step $1$: the construction of the AdS-spacetime with conformal 
boundary $(\widetilde{M}, \widetilde{O}, \widetilde{\mathcal M})$\/}

We need to start with a definition: an \emph{end\/} in $\widetilde{M}$ is an 
open domain $U$ in $\widetilde{M}$ such that the restriction of $\mathcal D$ 
to $U$ is injective, with image 
an end $V$ in $\overline{{\mathbb A}{\mathbb D}{\mathbb S}}$ 
(see definition~\ref{def.end}, Figure \ref{splitfigure}). 
The proof relies on the geometric understanding of ends, 
hence, we insist on their description: an end is the intersection 
$F \cap P$, where $P$ is the past of an element $a$ of 
$\overline{\mbox{Ein}}_2$, $F$ the future of an element $b$ 
of $\overline{\mbox{Ein}}_2$, and such that $a$, $b$ are strictly 
causally in a de Sitter domain containing $V$.
The interior in $\overline{\mbox{Ein}}_2$ of the intersection between $\overline{\mbox{Ein}}_2$ and the closure of $V$ is a diamond-shape region, denoted by $\partial{U}$. The end itself is the intersection between $\overline{{\mathbb A}{\mathbb D}{\mathbb S}}$ and the convex hull in $P(E)$ of the boundary of $\partial{U}$.

Let $\overline{U}$ be the union $U \cap \partial{U}$.
Observe that for any end $U$ the triple $(U, \partial{U},\overline{U})$ is 
an AdS-spacetime with boundary.
A \emph{marked\/} end is a pair $(U,x)$, where $U$ is an end in 
$\widetilde{M}$, and $x$ an element of 
$\partial {U} \subset \overline{\mbox{Ein}}_2$. 
Let $\mathcal E$ be the set of marked ends in $\widetilde{M}$. 
On $\mathcal E$, let $\sim$ be the equivalence relation identifying 
two marked ends $(U, x)$, $(U',x)$ if there is a third marked end $(U'',x)$ 
with $U'' \subset U \cap U'$. Let $\Upsilon$ be the quotient space of 
$\sim$. 
Let $\Xi$ be the union $\Upsilon \cup \widetilde{M}$.
For any end $U$, let $\partial\widehat{U} \subset \Upsilon$ be the set 
$\{ [U', y] \in \Upsilon / U' \subset U \}$, and let $\widehat{U}$ be 
the union in $\Xi$ of $\partial\widehat{U}$ with $U$. The 
$\partial\widehat{U}$ form the basis of a topology on $\Upsilon$, 
and the $\widehat{U}$ form, with the open subset of $\widetilde{M}$, 
the basis of a topology in $\Xi$.
It should be clear to the reader that in the special case where 
$\widetilde{M}$ is the end $V$,
all this process gives as final output topogical spaces $\Upsilon$, $\Xi$ respectively homeomorphic to $\partial{V}$, $V \cup  \partial{V}$, where $\partial{V}$ is the interior of the intersection between $\overline{\mbox{Ein}}_2$ and the closure of $V$.

The inclusion $\widetilde{M} \subset \Xi$ is clearly a homeomorphism onto 
its image, which is dense in $\Xi$. Similarly, any marked end $(U, x)$, 
the open domain $\widehat{U}$ is a neighborhood of $[U, x]$ in $\Xi$, 
which is homeomorphic to $V \cup \partial{V}$, the conformal completion 
in $\overline{\mbox{Ein}}_3$ of the end $V = {\mathcal D}(U)$ in 
$\overline{{\mathbb A}{\mathbb D}{\mathbb S}}$. It follows that 
$\Xi$ is a manifold, with chards the $\widehat{U}$ and the chards of 
$\widetilde{M}$. Indeed, the only remaining point to check 
(with the second countability that we leave to the reader) 
is the Hausdorff property: let $x_1$, $x_2$ be two elements of 
$\Xi$ such that every neighborhood of $x_1$ intersects every neighborhood 
of $x_2$. Then, clearly, if $x_1$ belongs to $\widetilde{M}$, 
the same is true for $x_2$, and $x_1 = x_2$. 
If not, $x_1$ and $x_2$ belong to $\Upsilon$: 
$x_1 = [ U_1, x'_1]$, $x_2 = [U_2, x'_2]$ with $(U_i, x_i) \in \mathcal E$. 
By hypothesis, the neighborhoods $\widehat{U}_1$ and $\widehat{U}_2$ must overlap. 
If $U_1 \cap U_2 \subset \widetilde{M}$ is empty, then some $[U_3, x_3]$ must belong 
to $\widehat{U}_1 \cap \widehat{U}_2$. Then, points in ${U}_3$ correspond to points 
in $\widehat{U}_3$ which are in $U_1 \cap U_2$: contradiction. Hence, $U_1 \cap U_2$ 
must intersect for any choice of the marked ends $(U_i, x_i)$. Fix one choice 
$(U^0_i, x^0_i)$ of these ends, and consider for every $i= 1,2$ smaller ends 
$U_i \subset U_i^0$. More precisely, fix $U_1$, with $[U^0_1, x_1] = [U_1, x_1]$ 
and $U_1 \subset U^0_1$. Then, if $U_2$ is sufficiently small, its image 
${\mathcal D}(U_2)$, which intersect ${\mathcal D}(U_1)$, is contained 
in ${\mathcal D}(U^0_1)$. According to Lemma~\ref{lem.fermassiette}, $U_2$ 
is then contained in $U^0_1$. But, since $x_1$, $x_2$ are not separated one 
from the other, they must have the same image under $\mathcal D$. 
Applying Lemma~\ref{lem.ferme}, we obtain $x=y$, i.e. $\Xi$ is Hausdorff.

The combination of the developing map 
${\mathcal D}: \widetilde{M} \rightarrow \overline{{\mathbb A}{\mathbb D}{\mathbb S}}$ 
with the inclusions $\overline{U} \subset \overline{\mbox{Ein}}_3$ 
induce a well-defined continuous map 
${\mathcal D}: {\Xi} \rightarrow \overline{\mbox{Ein}}_3$. 
Moreover, the restriction of ${\mathcal D}$ to (the projection of) 
any closed end $\overline{U}$ is injective: ${\mathcal D}$ is a local 
homeomorphism.

Finally, the action of $\Gamma$ on \mr extends naturally on ${\Xi}$: for any 
$\gamma$ in $\Gamma$, define  $\gamma[x \in \overline{U}]$ as being 
$[\rho(\gamma)x \in \rho(\gamma)\overline{U}]$.

Observe that this action is continuous and preserves 
$\Upsilon ={\Xi} \setminus \widetilde{M}$. Moreover, the map 
${\mathcal D}: {\Xi} \rightarrow \overline{\mbox{Ein}}_3$ is 
equivariant for this action.

We can now define $\widetilde{\mathcal M}$: this is the set of 
$\Gamma$-causally wandering points in ${\Xi}$, i.e. the set of elements 
$x$ of ${\Xi}$ admitting neighborhood $W$ such that for every non-trivial element 
$\gamma$ of $\Gamma$ no element of $W$ is causally related
in $\Xi$ to an element of $\gamma{W}$. Observe that ${\Xi}$ 
is open, 
$\Gamma$-invariant, and that it contains \mr since the action of 
$\Gamma$ on \mr is strongly causal. $\widetilde{O}$ is the 
complement $\widetilde{\mathcal M} \setminus \widetilde{M}$ 
(hence, the set of $\Gamma$-non causally wandering points).

\emph{Step $2$: the action of $\Gamma$ on $\widetilde{\mathcal M}$ is free, proper
and strongly causal.}

Observe that the action is free, since the fixed points 
are not wandering. Assume that the action of $\Gamma$ on $\widetilde{\mathcal M}$
is not proper.

Then, there are sequences $( x_n )_{(n \in {\mathbb N})}$, 
$( y_n )_{(n \in {\mathbb N})}$ in $\widetilde{\mathcal M}$ and a sequence $\gamma_n$ 
in $\Gamma$ such that:

-- $\gamma_nx_n = y_n$,

-- $x_n \to x \in \widetilde{\mathcal M}$,

-- $y_n \to y \in \widetilde{\mathcal M}$,

-- $y$ is not in the $\Gamma$-orbit of $x$.

Let $\bar{x} = {\mathcal D}(x)$, $\bar{y} = {\mathcal D}(y)$, $\bar{x}_n = {\mathcal D}(x_n)$, $\bar{y}_n = {\mathcal D}(y_n)$. We also decompose the $g_n = \rho(\gamma_n)$ along their left and right components: $g_n = (g_L^n, g_R^n)$.

Since the action of $\Gamma$ on \mr is proper (it is the group of covering automorphisms), $x$ or $y$ must belong to $\widetilde{O}$; let's say $x$. Then, $x$ belongs to
a diamond shape region $\partial U$, where $U$ is an end in $\widetilde{M}$. 
Since $x$ is $\Gamma$-wandering, we can choose $U$ so that 
$\gamma\overline{U} \cap \overline{U} = \emptyset$ for every non-trivial 
$\gamma$ in $\Gamma$. Moreover, since $\overline{U}$ is open in $\widetilde{\mathcal M}$, 
we can assume that the $x_n$ all belongs to $U$. 

Define $U_n = \gamma_nU$, $\Delta = \partial{U}$, 
$\Delta_n = \gamma_n\Delta$, and 
$\overline{\Delta}_n = {\mathcal D}(\Delta_n)$, $\overline{U}_n ={\mathcal D}(U_n)$.
The image $\overline{\Delta} = {\mathcal D}(\Delta)$ is a diamond-shape 
region. 
Since the $U_n$ are necessarily disjoint we can assume that none of them contain $y$.

\emph{Claim: $y$ belongs to $\widetilde{O}$.\/}

Assume not. There is a small neighborhood $W$ of $y$ in \mr 
such that the restriction of $\mathcal D$ is injective with image a small ellipsoid 
$\overline{W}$ in $\overline{{\mathbb A}{\mathbb D}{\mathbb S}}$. 
For $n$ sufficiently big $y_n$ belongs to $W$, hence
the intersection ${\mathcal I}_n = W \cap U_n$ is not empty. 
On the other hand, $\overline{U}_n$ is an end: it is a connected component 
of the intersection between $\overline{{\mathbb A}{\mathbb D}{\mathbb S}}$ 
and the complement of two hyperplanes in $P(E)$. It follows that the 
intersection $\overline{\mathcal I}_n$ between the ellipsoid $\overline{W}$ 
and $\overline{U}_n$ is convex: it is the trace in an ellipsoid of a half-space 
or a quarter of space. Moreover, the restriction of $\mathcal D$ to ${\mathcal I}_n$ and 
$W$ is injective. Hence, according to Lemma~\ref{lem.fermassiette}, 
the image of ${\mathcal I}_n$ by $\mathcal D$ is the entire 
$\overline{\mathcal I}_n$. In other words, the screen $\overline{W}$ 
reflects faithfully how the ends $U_n$ intersect $W$. But it is 
geometrically clear that $\overline{U}_n$ cannot accumulate to $\bar{y}$ 
if they are disjoint one to the other: visualize by considering an ellipsoid $W'$ in 
the Minkowski space conformally equivalent to $\overline{W}$; then, in this 
conformal chard, the $\overline{U}_n$ are intersections between 
the past of lightlike plane and the future of another lightlike plane. 
It leads to a contradiction: the claim is proved.

Replace the ellipsoid $W$ in the proof above by a neighborhood 
$\overline{U}'$ with $U' \in \mathcal E$ and such that $y$ belongs 
to $\partial{U}'$. We can assume that all the $y_n$ belong to $U'$. 
The argument above, based on lemma~\ref{lem.fermassiette}, 
shows that the image by $\mathcal D$ of the intersection $U_n \cap U'$ 
is also the entire $\overline{U}' \cap \overline{U}_n$. Applying once more 
this lemma to the closure, we obtain that $\mathcal D$ projects faithfully 
the intersections between $\Delta'$ and $\Delta_n$ over the entire 
$\overline{\Delta}_n \cap \overline{\Delta}'$. But these intersections 
are even simpler to visualize than intersection between an ellipsoid and 
ends: indeed, through the identification 
$\overline{{\mathbb A}{\mathbb D}{\mathbb S}} \approx {\mathbb R}P^1_L \times 
{\mathbb R}P^1_R$, a diamond-shape region corresponds to the product of two open 
intervals $I_L \times I_R$. Denote by $I^0_L$, $I^0_R$ the open intervals in
the projective line such that 
$I^0_L \times I^0_R = \partial\overline{\Delta}'$: $\bar{y}$ 
correspond to a pair $(y_L, y_R) \in I^0_L \times I^0_R$. 
For every integer $n$ let
$I^n_L \subset I^0_L$ and $I^n_L \subset I^0_R$ be the intervals such that
$\overline{\Delta}_n \cap \overline{\Delta}'= I^n_L \times I^n_R$.

Assume that
for some integers $n \neq m$ the intersection $I^n_L \cap I^m_L$ is not empty.
Then there is lightlike segment in $\overline{\Delta}'$ with one extremity
in $\overline{\Delta}_n$ and the other in $\overline{\Delta}_m$. It follows
that there is a causal curve joining an element of $U_n \cap U'$ 
to an element of $U_m\cap U'$. Now, for the first time, 
we use the fact that $M$ is strongly causal, i.e. that
the action of $\Gamma$ on $\widetilde{M}$ is strongly causal: it means
that the open domain $U$ can be selected so that for every non trivial
element of $\Gamma$ no element of $\gamma{U}$ can be causally related
to an element of $U$. Apply this remark to $\gamma_n^{-1}\gamma_m$: we
obtain a contradiction.

Hence, for every $n \neq m$ we have $I^n_L \cap I^m_L = \emptyset$.
Similarly $I^n_R \cap I^m_R =\emptyset$. But since the $y_n$ converge 
to $y = (y_L, y_R)$ it follows that $(I^n_L)_{(n \in {\mathbb N})}$ 
(resp. $(I^n_R){(n \in {\mathbb N})}$) is a sequence of intervals smaller
and smaller converging uniformly to $\{ y_L\}$ (resp. $\{ y_R \}$).
Hence for big $n$ we have $U_n \subset U'$. We obtain a contradiction
since $U'$ can be chosen so that $\gamma U'\cap U' =\emptyset$ for every non trivial
$\gamma$.

This final contradiction achieves the proof of step $2$. 
Hence, the quotient $\mathcal M =\Gamma\backslash\widetilde{\mathcal M}$ is a 
manifold, with boundary $O = \Gamma\backslash\widetilde{O}$.

\emph{Step $3$: $(M, O, {\mathcal M})$ is an universal conformal completion.}

Let $(M, O', {\mathcal M}')$ be another conformal completion of $M$. Let $\pi': \widetilde{\mathcal M}' \rightarrow {\mathcal M}'$ be the universal covering of 
$\mathcal M$.
The interior of $\widetilde{\mathcal M}'$ is simply-connected, hence it can be identified 
with the universal covering of $M$. Hence $M$, ${\mathcal M}'$ 
(and $\mathcal M$) have the same fundamental group $\Gamma$.
Let ${\mathcal D}': \widetilde{\mathcal M}' \rightarrow \overline{\mbox{Ein}}_3$ and 
$\rho: \Gamma \rightarrow G \times G$ be the developing map and the holonomy 
representation of $(M, O, {\mathcal M})$. 
The developing map for $M$ is then the restriction to $\widetilde{M}$ of 
${\mathcal D}'$, and $\rho$ is the holonomy representation of the AdS-spacetime $M$.

The main observation is the following: \emph{every point $x$ in $\widetilde{O}'$ 
admits a neighborhood $U'$ in $\widetilde{\mathcal M}'$ such that the restriction of 
${\mathcal D}'$ to $U'$ is injective and the image ${\mathcal D}'(U \cap \widetilde{M})$ 
is an end of AdS.\/} 
Hence, $(U \cap \widetilde{M}, {\mathcal D}'(x))$ is a marked end of $\widetilde{M}$. 
It defines a map $\widetilde{O} \rightarrow {\Upsilon}$. 
With the identity map $\widetilde{M} \rightarrow \widetilde{M}$, and after 
composition in the quotient space, we obtain a map 
$\widetilde{F}: \widetilde{\mathcal M}' \rightarrow {\Xi}$. 
The proof that $\widetilde{F}$ is a homeomorphism onto 
its image is straightforward and left to the reader. 

The main point is to show that $\widetilde{F}$ takes value in the domain of 
$\Gamma$-wandering points $\widetilde{\mathcal M}$. 
Assume by contradiction the existence of a point $x$ in $\widetilde{\mathcal M}'$ 
such that $x' = \widetilde{F}(x)$ is not a $\Gamma$-wandering point. 
Since $\Gamma$ acts properly on $\widetilde{\mathcal M}'$, 
there is a neighborhood $W$ of $x$ such that 
$\gamma{W} \cap W = \emptyset$ for every non-trivial $\gamma$. 
Moreover, since it is not $\Gamma$-wandering, $x'$ does not belong to 
$\widetilde{M}$. Hence, $x$ belongs to $\widetilde{O}'$: 
we can choose $W$ so that $W \cap \widetilde{M}$ is an end of 
$\widetilde{M}$.
Then, by hypothesis, there is a non-trivial $\gamma$ such that 
$\gamma\overline{W} \cap \overline{W}$ is not empty 
since $\overline{W}$ is a neighborhood of $x'$ in $\widetilde{\mathcal M}$. 
But completed ends overlapping in $\widetilde{\Xi}$ must also overlap in 
the interior of $\widetilde{\Xi}$, i.e. in $\widetilde{M}$. Hence, 
this overlapping exists also $\widetilde{\mathcal M}'$. Contradiction.

Hence, $\widetilde{F}$ takes value in $\widetilde{\mathcal  M}$. 
Since the construction is $\Gamma$-equivariant, $\widetilde{F}$ induces 
the required morphism $(M, {\mathcal M}', O') \rightarrow (M, O, {\mathcal M})$.
\fin

\rque
\label{rk.bordcausal}
Consider the quotient of $\widetilde{\mbox{AdS}}$ by a cyclic group $\Gamma$
generated by
a hyperbolic left translation $\gamma = (\gamma_L, id)$. Then the set of fixed points of 
$\gamma$ is the union of  two left lightlike leaves $l_1$, $l_2$. The action of $\Gamma$
on $\widetilde{\mbox{AdS}} \cup \widehat{\mbox{Ein}}_2 \setminus (l_1 \cup l_2)$ is free
and properly discontinuous, hence the quotient space is an universal conformal boundary
for $\Gamma\backslash\widetilde{\mbox{AdS}}$. But since every $\Gamma$-orbit of
elements of $\widehat{\mbox{Ein}}_2$ is contained in a right lightlike leaf, it follows that
the strongly causal conformal boundary of $\Gamma\backslash\widetilde{\mbox{AdS}}$ is empty.
\erque

\rque
\label{rk.nocompletion}
Restricting to the strongly causal category is essential to ensure the uniqueness of
the maximal conformal extension. Indeed: let $\gamma =
(\exp(\lambda\Delta),\exp(\mu\Delta))$ be a hyperbolic-hyperbolic element of
$\widetilde{G} \times \widetilde{G}$ with $\lambda \neq \mu$. The action of $\gamma$ on 
$\widetilde{\mbox{AdS}}$ is free and properly
discontinuous. In $\overline{\mbox{Ein}}_2$ t
here are two $\gamma$-invariant left leaves $l_1$, $l_2$, and
two $\gamma$-invariant right leaves $r_1$, $r_2$. 

The preimage $\tilde{l}_i$, $\tilde{r}_i$ in $\widehat{\mbox{Ein}}_2$ of  $l_i$, $r_i$
are lightlike geodesics. It is easy to see that the action of $\Gamma$ on 
$\widetilde{\mbox{AdS}} \cup \widehat{\mbox{Ein}}_2 \setminus (l_1 \cup l_2)$ is free
and properly discontinuous, and the same is true for the action on
$\widetilde{\mbox{AdS}} \cup \widehat{\mbox{Ein}}_2 \setminus (r_1 \cup r_2)$. 
Each of these extensions provides a conformal completion of
$\Gamma\backslash\widetilde{\mbox{AdS}}$. But it can be proved that each of this extension
is maximal, i.e. does not embeds in a larger conformal completion of 
$\Gamma\backslash\widetilde{\mbox{AdS}}$. Hence 
$\Gamma\backslash\widetilde{\mbox{AdS}}$ does not admit an universal 
conformal completion.

A similar situation appears in the so-called ``Taub - NUT'' examples, see \cite{taubnut, taubnut2}.

\erque

\section{BTZ black holes and multi black-holes}
\label{sec.last}

In this section, we consider various pairs $(\rho, {\Lambda})$, 
where $\rho = (\rho_L, \rho_R): \Gamma \rightarrow G \times G$ 
is an admissible representation and $\Lambda$ a $\rho$-admissible 
closed subset of $\overline{\mbox{Ein}}_2$ 
(see definition~\ref{def.admissible}).
We prove in the selected cases that the spacetimes 
$M_\Lambda(\Gamma)$ are AdS-spacetimes with black-holes 
in the meaning of definition~\ref{def.BH}.

\rque
\label{mollo}
We don't pretend to study all the possibilities: for example, for non-abelian $\Lambda$, we only consider the case $\Lambda = \Lambda(\Gamma)$, whereas the case $\Lambda \neq \Lambda(\Gamma)$ give other examples of spacetimes with black-holes.
Moreover, we don't consider the case where $\Gamma$ is trivial, which would lead, for the $\Lambda$ selected below, to AdS-spacetimes with black-hole too!

A justification for this omission is that all these spacetimes are not 
maximal for (too much) obvious reasons, because of the embedding 
$M_\Lambda(\Gamma) \subset M_{\Lambda(\Gamma)}(\Gamma)$ - especially in the case $\Gamma = id$, where the whole spacetime embeds in $\widetilde{\mbox{AdS}}$. 
\erque

\subsection{The conical black-holes}

This is the case where $\Lambda$ is conical (see \S~\ref{subsub.conical}). More precisely,
we have to consider the case where $\Lambda$ is a upper lower tent ${\mathcal T}^+_{xy} = [x,z] \cup [z,y]$ (if it is a lower tent, then the spacetime is full, but without black-hole).
The domain $E(\Lambda)$ has been described in \S~\ref{subsub.conical}. From this description (see also 
the Figure \ref{splitfigure}), it appears clearly that the conformal boundary $\widetilde{O}$ of $E(\Lambda)$ is the diamond-shape region $\Delta_2 = \Omega(\Lambda)$, and the full completion of $E(\Lambda)$ is ${\mathcal M}(\Lambda) = E(\Lambda) \cup \Delta_2$. The region $F_1 \cap F_2$ is the region invisible from $\widetilde{O}$. We can understand a part of the terminology reported in remark~\ref{rk.BTZnotations}: the outer region $F_2 \cap P_2$ is the region visible by observers in $O$.

According to \S~\ref{conicalbof} the action of $\Gamma$ on $E(\Lambda)$ is free, properly discontinuous and strongly causal if and only if $\Gamma$ is a cyclic group generated by an element $\gamma = (\gamma_L, \gamma_R) = (\exp(u\Delta), \exp(v\Delta))$. We can assume without loss of generality $v \geq u \geq 0$. It appears then clearly that the action of $\Gamma$ on $\widetilde{O} = \Delta_2$ is free and proper. Moreover, this action is strongly
causal if $u=0$.
Hence, its quotient $O$ is the natural conformal boundary of $M_{\Lambda}(\Gamma)$
and it is also the strongly causal conformal boundary if $u \neq 0$ - if $u=0$, the strongly
causal conformal boundary is empty (see remark~\ref{rk.bordcausal}).
Obviously, the invisible domain in $M_\Lambda(\Gamma)$ from $O$ is the quotient of the ``intermediate region'' $F_1 \cap F_2$. 

It follows that $M_\Lambda(\Gamma)$ is an AdS-spacetime with black-hole in
the meaning of definition~\ref{def.BH} if and only if $u \neq 0$. 

\rque
\label{rk.topoconical}
The topology is very simple: $M_\Lambda$ is homeomorphic to the product of the annulus by $\mathbb R$. The same is true for the outer region and the black-hole, which are separated 
by a lightlike annulus, the horizon. 
\erque

\subsection{Splitting black holes:}

Here, we consider the case where $\Lambda$ is splitting, i.e. 
two non-causally related points $(x,y)$. This case is fully detailed 
and described in \S~\ref{subsub.split}. Figure \ref{splitfigure} 
is still useful.
The discussion above remains essentially the same, but now the 
conformal boundary of $E(\Lambda)$ has two connected components: 
$\Delta_1$ and $\Delta_2$. Hence, the invisible domain from their 
union is still $F_1 \cap P_1$.

According to \S~\ref{splitcausal}, 
in order to act properly on $E(\Lambda)$, 
the group $\Gamma$ must be as in the conical case: generated by a hyperbolic translation, or a hyperbolic-hyperbolic element. Hence, the conformal completion of the quotient is the quotient of these two diamond-shapes regions. 
The strongly causal conformal boundary is $\Delta_1 \cup\Delta_2$, except if
$\gamma$ is an hyperbolic translation, in which case the strongly causal boundary is
empty.

When $\gamma$ is hyperbolic-hyperbolic, 
$M_{\Lambda}(\Gamma)$ contains a black-hole - the quotient of $F_1 \cap F_2$ - 
isometric to the black-hole of the conical case. 

The topological description is the same than in the conical case. But the horizon is  not
$C^2$, and the visible domain is not globally hyperbolic.

\rque
\label{rk.singleBTZ}
Something similar to what was discussed in remark~\ref{mollo} appears: 
the conical spacetime $M_\Lambda(\Gamma)$ embeds isometrically in 
$M_{xy}(\Gamma)$: hence, it is not maximal. 
This is a good reason for considering that conical 
case does not contain a black-hole, as in the ``classical'' 
litterature. But observe that $M_{xy}(\Gamma)$ itself is not maximal 
too: when $u < v$, $M_{xy}(\Gamma)$ 
embeds in the spacetime $M^D(\Gamma) = \Gamma\backslash{D}(\Gamma)$, 
where $D(\Gamma)$ is the absolute 
causality domain (see definition~\ref{def.absolu}, 
proposition~\ref{prop.Dcausal}, 
case $(6)$). Observe that $D(\Gamma)$ is an open domain in 
$\widetilde{\mbox{AdS}}$ 
not contained in a affine domain. Actually, the proof of case $(6)$ of Proposition~\ref{prop.Dcausal} shows that $D(\Gamma)$ contains all the 
preimage in $\widetilde{\mbox{AdS}}$ of $E_{xy}$, which are connected 
one to the other by regions (the ``inner regions'' with the conventions in 
\cite{BTZ, BTZ2, brill}) which project in 
$\mbox{AdS} \approx \mbox{SL}(2,{\mathbb R})$ to the domain 
$0 < bc < \frac{(\exp(u)-\exp(v))^2}{4\exp(u+v)}$. Observe that these 
new regions has empty conformal boundary: their closure intersect 
$\mbox{Ein}_2$ only along the union of two lightlike geodesics. 
Hence, the conformal boundary of $D(\Gamma)$ is the preimage in 
$\widehat{\mbox{Ein}}_2$ of all the preimage of $\Delta_1 \cup \Delta_2$. 

In other words, let $A$ is an affine domain containing $E(\Lambda)$, 
and let $A_i$ be the infinite family of affine domains in 
$\widetilde{\mbox{AdS}}$ such that $\widetilde{\mbox{AdS}}$ 
is the union of the $\overline{A}_i$ (see \S~$3.5$ in \cite{ba1}). 
Every $\overline{A}_i$ contains a copy $E_i$ of $E(\Lambda)$.
Moreover, $A$ can be selected so that the conformal boundary of $\overline{A}_i$ 
contains the conformal boundary ${\Delta}_1^i \cup \Delta_2^i$ of $E_i$. 
The inner regions connect every $\overline{A}_i$ to the following 
$\overline{A}_{i+1}$, offering the way to some causal curves to pass 
from $\overline{A}_i$ to $\overline{A}_{i+1}$.
It is easy to show that, thanks to these connecting inner regions, $D(\Gamma)$ 
is entirely visible from its conformal boundary: it has no hole. 

However, the quotient $M^D(\Gamma)$ is considered in the litterature devoted 
to BTZ black-holes (including \cite{BTZ}) has the typical 
spacetime containing a \emph{single non-static ($u \neq v$) BTZ-black-hole.\/} 
It means that the point of view to adopt is to pay attention to simple blocks 
$\overline{A}_i \cap D(\Gamma)$ individually, to consider the 
observers only in the boundary components of \emph{one} of them, 
and to consider other blocks as being other parts of the universe 
which can be reached only by going through the horizon of the black-hole.
Adopting this point of view, we observe that $M_{xy}(\Gamma)$ is 
therefore enough to give a picture of the considered black-hole.

\rque
This causal description is actually very similar to the description of black-holes in the maximal Kerr spacetime $M^{max}_{Kerr}$ (see the Introduction).
\erque

\erque

\subsection{The extreme black hole}
\label{sub.extremebtz}

The case where $\Lambda$ is extreme is described in \S~\ref{subsub.extreme}. 
In this case, the invisible domain is not contained in an affine domain, hence we need
two successive domains, and $\Lambda$ must be considered as a closed subset 
of $\widehat{\mbox{Ein}}_2$ (see Figure \ref{extremefigure}). In this case, 
the conformal boundary is the extreme diamond $\Omega(x,y)$, 
and the entire $E(\Lambda)$ is visible from the boundary:  
there is no black-hole!

According to \S~\ref{sub.invextreme}, the group $\Gamma$ must 
be contained in $A_{hyp}$ or $A_{ext}$, and since we want the action of 
$\Gamma$ 
on $E(\Lambda)$ to be causal, the case $\Gamma \subset A_{hyp}$ must be excluded (\S~\ref{subsub.extremehyp}).
Hence, the action of $\Gamma$ on $E(\Lambda)$ is free, properly discontinuous and 
strongly causal, except if 
it contains a parabolic translation (\S~\ref{subsub.extremepar}).

Anyway, as we have seen, $E(\Lambda)$ does not contain any black-hole. 
But we can use the
same trick as for single non-static BTZ black-holes (remark~\ref{rk.singleBTZ}): 
consider 
the absolute causal domain $D(\Gamma)$. Observe that $D(\Gamma) = E(\Lambda)$ if $\Gamma$ 
is not cyclic (Theorem~\ref{thm.egal}). Hence, this trick will apply only for cyclic subgroups
of $A_{ext}$. Moreover, according to Proposition~\ref{prop.Dcausal}, case $(3)$, the 
absolute causal domain of parabolic translations is trivial: we must exclude them. Finally,
elements of $\Gamma$ are not hyperbolic translations since $\Gamma \subset A_{hyp}$ is excluded.
Hence, $\Gamma$ must be generated by a parabolic-hyperbolic element. According to 
Proposition~\ref{prop.Dcausal}, case $(7)$, the action of $\Gamma$ on $D(\Gamma)$ is
free, properly discontinuous and proper. 

We can reproduce nearly the same comment than in remark~\ref{rk.singleBTZ}: $D(\Gamma)$ has to
be understood as a $\delta$-invariant subset of $\widetilde{\mbox{AdS}}$, and the quotient 
space $M^D(\Gamma)$ is an union of ``local universes''. There is a small difference: 
the conformal boundary of every simple block is now connected 
(the other connect component ``vanished''), and there is no intermediate region: the black-holes
correspond to inner regions.

\subsection{Multi-black holes}
The last case we consider is the non-elementary case: $\Gamma$ is non-abelian, and 
$\Lambda = \Lambda(\Gamma)$ is not elementary. The key point is to use 
Proposition~$8.50$ of \cite{ba1}: \emph{$E(\Lambda)$ is the union of 
the past and future globally hyperbolic cores, with the closed ends $\Omega(I)$ associated to 
gaps $I$ of $\Lambda$.\/} Observe that the gaps in this case are not extreme 
(remark~$8.25$ of \cite{ba1}),
hence, the diamonds $\Omega(I)$ are not extreme.
Moreover, since in this case the left and right morphisms are both faithfull,
the stabilizer of $\Omega(I)$ is generated by a hyperbolic - hyperbolic element.

It follows clearly that the connected components of the conformal boundary 
of $E(\Lambda)$ 
are precisely the diamonds $\Omega(I)$: the only way to approach this 
conformal boundary is to
enter in a closed end. Observe that the closed end associated to a gap $I$ is isometric
to the outer region of a conical spacetime. Moreover, if $\Lambda$ is a topological circle, 
then the conformal boundary is empty: there is no observor, no black-hole.
Finally, the invisible domain from $\Omega(\Lambda)$ is precisely the future 
globally hyperbolic core $E(\Lambda^+)$.

Now, according to Theorem~\ref{noneleOK}, and since every element of $\Gamma$ 
is hyperbolic - hyperbolic, the action of $\Gamma$ on 
$\Omega(\Lambda) \cup E(\Lambda)$
is free, properly discontinuous and strongly causal 
(observe that it is true even if $\Lambda \neq \Lambda(\Gamma)$). 
The quotient 
spacetime $M_\Lambda(\Gamma)$ is then a AdS-spacetime with one black-hole: the 
(globally hyperbolic) quotient of $E(\Lambda^+)$ by $\Gamma$.

\rque
There is an obvious $1$-$1$ correspondance between the connected 
components of the conformal 
boundary of $M_\Lambda(\Gamma)$ and $\Gamma$-orbits 
of gaps, i.e. non-cuspidal boundary components of 
the surfaces $\Sigma_L \approx \Sigma_R$ (see remark~\ref{rk.endpara}).
\erque

\rque
According to Theorem~\ref{thm.egal}, the trick used in 
remark~\ref{rk.singleBTZ} in order 
to enlarge the spacetime by considering the absolute causal 
domain $D(\Gamma)$ gives nothing new:
the quotient is the same spacetime.
\erque

\rque
\label{rk.pasmaximal}
Anyway, even if the ``obvious'' trick above does not work, many spacetimes $M_{\Lambda}(\Gamma)$ are
not maximal AdS-spacetimes. Indeed, add to $E(\Lambda)$ a 
very small end $U$,
such that the intersection $\partial U$ of the closure of $U$ with 
$\overline{\mbox{Ein}}_2 \approx {\mathbb R}P^1_L \times {\mathbb R}P^1_L$ is a 
small
rectangle $I_L \times I_L$ around a point $x$ in the boundary of $\Omega(I)$ (for some gap $I$).
Don't take as point $x$ one of the corner points of $\Omega(I)$. Then, if for some element 
$\gamma$ of $\Gamma$ the intersection $\partial{U} \cap \gamma\partial{U}$ is not empty, then
$\gamma\Omega(I) \cap \Omega(I)$ is not empty. Then, $\gamma$ must belong to the cyclic
subgroup $\Gamma_I$ of elements preserving the gap $I$ (see the proof of Theorem~\ref{noneleOK}).
Thus, $U$ can be chosen so that the intersection $\gamma{U} \cap U$ never happens. Then,
the union of $E(\Lambda)$ with all the $\gamma{U} \;\; (\gamma \in \Gamma)$ is 
a $\Gamma$-invariant connected 
spacetime $E'$ on which $\Gamma$ acts freely and properly discontinuously. 
Its quotient is a bigger 
spacetime $M'$ in which $M_\Lambda(\Gamma)$ embeds isometrically. 
\erque

\rque
\label{rk.excausal}
Consider once more the spacetime $M'$ constructed in the remark above.
Let's prove that $M'$ is strongly causal.
There is no loss of generality in
assuming that the point $x$ has been selected so that its (local) future does
not meet $E(\Lambda) \cup \Omega(\Lambda)$. Then, no future oriented causal
curve starting from a point in $E' \setminus E(\Lambda)$ can enter into $E(\Lambda)$:
we say that $E(\Lambda)$ is \emph{past-convex\/} in $E'$.

Let $y$ be an element of $E'$. If $y$ belongs to $E(\Lambda)$, since $E(\Lambda)$
is causally convex and strongly causal, there is a neighborhood $V$ of $y$ which is 
not causally related to any non-trivial $\Gamma$-iterate of itself.
Assume now that $y$ belongs to the end $U$ but not to $E(\Lambda)$. 
Observe that by construction $U$ meets $E(\Lambda)$ only in the closed end 
$E(I)$ associated to $I$. It follows that small neighborhoods of $y$ are contained
in $U \setminus E(\Lambda)$. Any causal curve $c$ joining a point near $y$ to a point 
near $\gamma y$ with $\gamma \neq id$ must enter in $E(\Lambda)$.
Since $E(\Lambda)$ is past-convex in $E'$ the causal curve $c$ 
must be past-oriented. 
But the same argument applied
near $\gamma y$ shows that $c$ must be future oriented.

This contradiction implies that the action of $\Gamma$ on $E'$ 
 is strongly causal. Therefore, $M'$ is strongly causal.

It is quite clear that the conformal boundary of $M'$ is the quotient by $\Gamma$
of the union of $\Omega(\Lambda)$ with the disjoint union of all the 
$\gamma\partial{U}$. Moreover, by causal convexity of 
$E(\Lambda) \cup \Omega(\Lambda)$, the strongly causal conformal boundary 
of $M'$ contains as an open subset
the conformal boundary of $M_\Lambda(\Gamma)$.
We claim that the conformal boundary of $M_\Lambda(\Gamma)$
is actually a connected component of the strongly causal conformal 
boundary of $M'$. Indeed, for every $y$ in $\partial{U}$ and in the boundary
of $\Omega(I)$ in the Einstein space, for 
any neighborhood $V$ of $y$ in $\partial{U}$ and for any non-trivial element
$\gamma$ of $\Gamma_I$ there are elements of $V \cap \Omega(I)$ causally 
related to elements of $\gamma{V} \cap \partial{U}$. 

It follows quite easily that the strongly causal conformal boundary of 
$M'$ is the quotient by $\Gamma$ of the union of $\Omega(\Lambda)$ with the
disjoint union of all $\gamma (\partial{U} \setminus \overline{\Omega}(I))$, where
$\gamma$ describes all $\Gamma$, and 
$\overline{\Omega}(I)$ is the closure in the Einstein space of $\Omega(I)$.

Points in $E' \setminus E(\Lambda)$, i.e. elements of $\gamma{U}
\;\; (\gamma \in \Gamma)$ are in the past of the strongly causal boundary of
$E'$. Hence, the black-hole in $M'$ is contained in the black-hole of 
$M_\Lambda(\Gamma)$, in particular, smaller. Moreover, 
if the point $x$ has been selected in the past of $\Omega(I)$ (the proof of
strong causality still applies, just replace past-convex above 
by future-convex)
then the black-hole in $M'$ is equal to the black-hole in $M_{\Lambda}(\Gamma)$:
the part of the conformal boundary we added didn't reveal any new point.
But if the point $x$ has been selected in the future $\Omega(I)$ there is no 
general answer: the black-hole in $M'$ might be strictly smaller than
the black-hole in $M_\Lambda(\Gamma)$.
\erque

\rque
\label{rk.bordcausalgh}
The example described in remarks~\ref{rk.pasmaximal}, \ref{rk.excausal}
shows that spacetimes with BTZ multi black-holes are far from being 
maximal as spacetimes, even in the strongly causal category. But observe that these
examples become forbidden if in the definition of
spacetimes with black-holes we impose the following additionnal requirement:
\emph{the strongly causal conformal boundary must be globally hyperbolic.\/}
\erque

\subsection{Kerr-like coordinates}
\label{sub.Kerrlike}

Contrary to the habits in the BTZ litterature, we end by the presentation of 
the Kerr-like expression of the BTZ metric. 
This expression concerns the metric on the ``outer regions'', i.e. the ends. 
Hence we just have to consider the splitting and extreme cases.

\subsubsection{The Kerr-like metric on the outer region of a splitting spacetime.}

In this case, $\Gamma$ is generated by $\gamma$, where, up to conjugacy:

\[ \gamma = (\gamma_L, \gamma_R) = (\exp(u\Delta), \exp(v\Delta)), \;\; v \geq  u > 0 \]
The elements $x$, $y$ of $\Lambda$ 
must be the attractive and repulsive fixed points of $\gamma$.

Define $r_\pm$ ($r_+> r_- \geq 0$) by:

\[ u = \pi(r_+ - r_-) \;\;\; v = \pi(r_+ + r_-) \]

Then: 

\[ \gamma_L g \gamma_R^{-1} = \left(\begin{array}{cc}
               a\exp(-2\pi r_-) & b\exp(2\pi r_+) \\
               c\exp(-2\pi r_+) & d\exp(2\pi r_-) \end{array}\right) \]

Let $(x_1, x_2, y_1, y_2)$ be coordinates of $E$ for which $Q= dx_1^2 + dx_2^2 - dy_1^2 - dy_2^2$.
Consider the identification:

\begin{eqnarray*}
\mbox{AdS} & \rightarrow & \mbox{SL}(2, {\mathbb R}) \\
(x_1, x_2, y_1, y_2) & \mapsto & \left(\begin{array}{cc}
               y_1+x_1 & x_2+y_2 \\
               x_2-y_2 & y_1-x_1 \end{array}\right)
\end{eqnarray*}

Then, in the coordinates $(x_1, y_1, x_2, y_2)$, the matrix in  $\mbox{SO}(2,2)$ corresponding to $\gamma$ is:

\[ \left(\begin{array}{cccc}
 \cosh(2\pi r_-) & -\sinh(2\pi r_-) & 0 & 0 \\
-\sinh(2\pi r_-) & \cosh(2\pi r_-) & 0 & 0 \\
0 & 0 & \cosh(2\pi r_+) & \sinh(2\pi r_+) \\
0 & 0 & \sinh(2\pi r_+) & \cosh(2\pi r_+)
\end{array}\right) \]

The attractive and repulsive fixed points have coordinates $(0,\pm 1,0 ,1 )$.
Thus, the outer region is $\{ \mid y_1 \mid  < x_1, \mid x_2 \mid < y_2 \}$. 
A natural associated coordinate system on this domain is:

\begin{eqnarray*}
x_1 & = & \rho_1 \cosh(T) \\
y_1 & = & \rho_1 \sinh(T) \\
x_2 & = & \rho_2 \sinh(\phi) \\
y_2 & = & \rho_2 \cosh(\phi)
\end{eqnarray*}

with $\rho_1^2 = \rho_2^2 - 1$.
The AdS-metric $dx_1^2 + dx_2^2 - dy^2_1 - dy_1^2$ in the coordinates $(T,\phi, \rho = \rho_2)$ is:

\[ \frac{1}{\rho^2 - 1}d\rho^2 + \rho^2d\phi^2 - (\rho^2-1)dT^2  \]

Observe that $\phi$, $T$ may have any real value, and that $\rho$ takes value in $]1,+\infty[$.

The action of $\gamma$ in the coordinates $(\phi, T, \rho_1, \rho_2)$ 
is simply the translation by $2\pi r_+$ on $\phi$ and by 
$-2\pi r_-$ on $T$. Hence, it preserves the function 
$t = \frac{r_-\phi + r_+T}{r_+^2 - r_-^2}$, and adds to 
$\varphi =\frac{r_+\phi + r_-T}{r_+^2 - r_-^2}$ the term $2\pi$. 
Therefore, we introduce the coordinates $(t, \varphi)$ instead of $(T, \phi)$:
$\varphi$ is considered as a polar coordinate, the action by $\gamma$ being 
represented by $\varphi \to \varphi+2\pi$, the other coordinates 
$(t, \rho)$ remaining unchanged. Actually, replace $\rho$ 
by the coordinate $r$ defined by: $\rho^2 = \frac{r^2 - r_-^2}{r_+^2 - r_-^2}$. 
Then, the AdS-metric becomes:

\[ -N(r)dt^2 + N(r)^{-1}dr^2 + r^2(d\varphi + \frac{J}{2r^2}dt)^2 \]

where:
\begin{eqnarray*}
J & = & -2r_-r_+ \\
N(r) & = & \frac{(r^2 - r_+^2)(r^2 - r_-^2)}{r^2} 
\end{eqnarray*}

The coordinates $(\varphi, t, r)$ are the \emph{Kerr-like coordinates.\/} Considering $\varphi$ as defined modulo $2\pi$, it provides an expression of the outer region of $M_{xy}(\Gamma)$.
The analogy with the Kerr metric is striking if we observe $N(r) = r^2 - M + (\frac{J}{2r})^2$, where $M = r^2_++ r_-^2$.

\rque
\label{rk.maximal}
In the Kerr-like coordinates, the level set of the time function $t$ are \emph{not} homogeneous.
The stabilizer of the outer region is $A_{hyp}$, the orbits of which in the outer region are timelike. Observe that this action of $A_{hyp}$ is an action by translations on the coordinates $(t, \varphi)$. Actually, the (scalar) curvature of the level sets of $t$ is not constant (except if $r_-=0$). What is the specific geometric feature of these level sets?

A remarkable fact is that the level sets $\{ t=t_0\}$ are \emph{maximal,\/} i.e. have zero mean curvature. A quick way to perform the computation is as follows: parametrize the level set $S_0 = \{ t=t_0 \}$ by parameters $\eta, \varphi$, where $\rho = \cosh(\eta)$:

\begin{eqnarray*}
x_1 & = & \sinh(\eta)\cosh(T) \\
y_1 & = & \sinh(\eta)\sinh(T) \\
x_2 & = & \cosh(\eta)\sinh(\phi) \\
y_2 & = & \cosh(\eta)\cosh(\phi)
\end{eqnarray*}

with $T = r_+t_0 - r_-\varphi$, $\phi = r_+\varphi - r_-t_0$.
Let $p$ be a point of $S_0$ of coordinates $(\eta, \varphi)$. Identify the tangent space to AdS at $p$ with the $Q$-orthogonal $p^\perp$ in $E$. Then, the tangent vectors of $S_0$ at $p$ are generated by:

\[ 
\partial_\eta = \left(\begin{array}{c}
\cosh(\eta)\cosh(T) \\
\cosh(\eta)\sinh(T) \\
\sinh(\eta)\sinh(\phi) \\
\sinh(\eta)\cosh(\phi)
\end{array}\right), \;\;\;
\partial_\varphi = \left(\begin{array}{c}
-r_-\sinh(\eta)\sinh(T) \\
-r_-\sinh(\eta)\cosh(T) \\
r_+\cosh(\eta)\cosh(\phi) \\
r_+\cosh(\eta)\sinh(\phi)
\end{array}\right)
\]

The future oriented normal $n_0$ to $S_0$ at $p$ is the following vector $n$, divided by the square root of the opposite of its norm, which is $r_-^2\sinh(\eta)^2 - r_+^2\cosh(\eta)^2 = -r^2$:

\[n = \left(\begin{array}{c}
r_+\cosh(\eta)\sinh(T) \\
r_+\cosh(\eta)\cosh(T) \\
-r_-\sinh(\eta)\cosh(\phi) \\
-r_-\sinh(\eta)\sinh(\phi)
\end{array}\right)
\]

Now, the second fundamental form of $S_0$ on a tangent vector field $X$ at $p$ is obtained by computing $\langle n_0 \mid \nabla_XX \rangle$. But the Levi-Civita connection of AdS is just the orthogonal projection on TAdS of the usual flat connection of $E$. Hence, the (extrinsic) curvature of the curves $\{ \varphi= Cte\}$ and $\{ \eta = Cte \}$ are the $Q$-scalar products with $n_0$ of the following vectors:

\[ 
\partial_{\eta\eta} = \left(\begin{array}{c}
\sinh(\eta)\cosh(T) \\
\sinh(\eta)\sinh(T) \\
\cosh(\eta)\sinh(\phi) \\
\cosh(\eta)\cosh(\phi)
\end{array}\right) = p,\;\;\;
\partial_{\varphi\varphi} = \left(\begin{array}{c}
r_-^2\sinh(\eta)\cosh(T) \\
r_-^2\sinh(\eta)\sinh(T) \\
r^2_+\cosh(\eta)\sinh(\phi) \\
r_+^2\cosh(\eta)\cosh(\phi)
\end{array}\right)
\]
These scalar products are null. Since the curves $\{ \eta  = Cte \}$ and $\{ \varphi=Cte\}$ are everywhere orthogonal, the mean curvature of $S_0$ at $p$ is thus $0$ (but these curves do not define the principal directions if $r_- \neq 0$). Observe that $\{ \varphi = Cte \}$ is actually a spacelike geodesic in AdS.

Moreover, $\langle n \mid \partial_{\eta\varphi} \rangle = r_-r_+$. It follows that the second fundamental form is $\mbox{II} = r^{-1}r_-r_+d\eta{d}\varphi$. The pair 
$(\partial_\eta, r^{-1}\partial_\varphi)$ is an orthonormal basis of the tangent space at $p$. Hence, the Gauss curvature of $S_0$ is:

\[ \frac{r_-^2r_+^2}{4r^4} \]

\erque

\subsubsection{The Kerr-like metric on the outer region of an extreme spacetime.}

Observe that the Kerr-like metric remains meaningfull when $r_-= r_+$, even if the 
coordinate transformations considered in the previous {\S} are not anymore valid.
Let $E'(r_+)$ be (simply connected) lorentzian manifold
consisting of ${\mathbb R}^2 \times ]r_+, +\infty[$ with coordinates $(t, \varphi, \rho)$,
equipped with the metric:

\[ -N(r)dt^2 + N(r)^{-1}dr^2 + r^2(d\varphi + \frac{J}{2r^2}dt)^2 \]

where:
\begin{eqnarray*}
J & = & -2r_+^2 \\
N(r) & = & \frac{(r^2 - r_+^2)^2}{r^2} = r^2 - 2r_+^2 + (\frac{J}{2r})^2
\end{eqnarray*}

$M'(r_+)$ is the quotient of $E'(r_+)$ by the translation $\varphi \to \varphi+2\pi$, the other
coordinates remaining the same. We want to prove that $E'(r_+)$ is isometric to the outer region
of an extreme black-hole as defined in {\S}~\ref{sub.extremebtz}.

The sectional curvature of the Kerr-like metric for $r_- = r_+ -\epsilon$ is $-1$ for any $\epsilon$, it then remains true at the limit $\epsilon = 0$. Hence, $E'(r_+)$ and $M'(r_+)$ 
are locally AdS. Observe also that $\partial_\varphi$ and $\partial_t$ generates a rank $2$ abelian Lie algebra of Killing vector fields. This is the pull-back by the developing map of an
abelian Lie algebra ${\mathcal A}(r_+)$ of Killing vector fields on AdS. This Lie algebra is
of course a limit of algebras ${\mathcal A}(r_+, r_+ -\epsilon)$ which are all Lie algebras of
subgroups conjugate to $A_{hyp}$. It follows easily that ${\mathcal A}(r_+)$ is either conjugate to the Lie algebra of $A_{hyp}$, or to the Lie algebra of $A_{ext}$. 
A quick calculus shows that $\partial_\varphi + \partial_t$ is an everywhere lightlike Killing vector field: it excludes the hyperbolic-hyperbolic case $A_{hyp}$, hence, up to conjugacy, the isometry group generated by $\partial_t$ and $\partial_\varphi$ is $A_{ext}$. 

Parametrize a line $\{ t=Cte, \rho=Cte \}$ by $\varphi$, and compute $\langle \partial_{\varphi\varphi} \mid p \rangle = -r^2$  
(once more, we can compute for $r_- = r_+ - \epsilon$). Hence, the orbits of the translations on the $\varphi$-coordinate are not geodesic. It follows that they are not hyperbolic translations: they are hyperbolic-parabolic. In other words, the monodromy of the translation $\varphi \to \varphi+2\pi$ is conjugate to $\gamma' = (\exp(2{\pi}r_0\Delta),\exp(H))$, for some $r_0 > 0$.
Actually, $r_0 = r_+$: indeed, for $r_- = r_+ - \epsilon$, the left component of the monodromy
of $\varphi \to \varphi +2\pi$ is conjugate to $\exp(2{\pi}r_+)$, it remains true at the limit $\epsilon = 0$ by continuity of the monodromy under deformations of the AdS-structure.
Hence, $\gamma' = \gamma$.

We observe now that $\partial_\varphi$ is spacelike on $E'(r_+)$. Hence, the image of the developing map of $E'(r_+)$ is contained in $D(\gamma)$. The curve $c(r) = (t=0, \varphi=0, r)$  is a geodesic (for example, the study in the preceding section shows that it is true in the case $r_- < r_+$, and the case $r_-=r_+$ is a limit case).
Moreover, the length between two points $c(r)$, $c(r')$ is $\frac{1}{2}\log(\frac{r'^2-r_+^2}{r^2-r_+^2})$: $c$ is a complete spacelike geodesic. According to the description
of absolute causal domains of hyperbolic-parabolic elements (case $(7)$ of Proposition~\ref{prop.Dcausal}), the complete spacelike geodesics entirely contained in
$D(\gamma)$ must actually lie in $C(\gamma)$, which is the outer region $E$ of the extreme black-hole associated to $\gamma$. Moreover, every $A_{ext}$-orbit in $E'(r_+)$ intersects $c$. Hence, the developing image of $E'(r_+)$ is contained in $E$. The action of $A_{ext}$ on
$E$ is free: it follows that the restriction of the developing map $\mathcal D$ to any $A_{ext}$-orbit in $E'(r_+)$ is a homeomorphism onto an entire $A_{ext}$-orbit in $E$.
Finally, the restriction of  $c$ is injective, with image a complete spacelike geodesic $\bar{c}$ in $E$. Every $A_{ext}$-orbit in $E$ intersect $\bar{c}$. It follows that $\mathcal D$ is an isometry between $E'(r_+)$ and $E$: the coordinates $(t, \varphi, r)$ parametrizes the entire
outer region $E$. 

\rque
The limit ``$r=r_+$'' is a lightlike segment $\tau$ in $\partial\mbox{AdS}$. It follows that the time levels $\{ t=Cte \}$ are closed in AdS. Their closure in $\mbox{AdS} \cup \partial\mbox{AdS}$ intersects $\partial\mbox{AdS}$ along an achronal topological circle which contains the lightlike segment $\tau$, and is spacelike outside $\tau$. Compare this situation with the example described page $45$ in \cite{mess}.

Of course, these level sets have zero mean curvature since it is true in the non-extreme case $r_- < r_+$. Furthermore, they have Gauss curvature $r_+^4/4r^4$.
\erque


\begin{thebibliography}{100}

\bibitem{amaking} 
S. Aminneborg, I. Bengtsson, S. Holst, 
P. Peldan, \emph{Making anti-de Sitter Black holes,\/}  
Class. Quant. Grav. , {\bf 13}  (1996), no. 10, 2707--2714, gr-qc/9604005.

\bibitem{aminneborg} 
S. Aminneborg, I. Bengtsson, D.R. Brill, S. Holst, 
P. Peldan, \emph{Black holes and wormholes in $2+1$ dimensions,\/}  
Class. Quant. Grav. , {\bf 15}  (1998),  no. 3, 627--644,
gr-qc/9707036.

\bibitem{amispin} 
S. Aminneborg, I. Bengtsson, S. Holst, \emph{A spinning anti-de 
Sitter wormhole,\/} Class. Quant. Grav., {\bf 16}  (1999), no. 2, 363--382, gr-qc/9805028 v2.


\bibitem{ashtekar} 
A. Ashtekar, G.J. Galloway, \emph{Some uniqueness results for dynamical horizons,\/} gr-qc/0503109


\bibitem{BTZ} 
M. Ba\~{n}ados, C. Teitelboim, J. Zanelli, \emph{The Black hole in
three-dimensional spacetime,\/} Phys. Rev. Lett., {\bf 69} no 13 (1992),
1849--1851, hep-th/9204099.

\bibitem{BTZ2} 
M. Ba\~{n}ados, M. Henneaux, C. Teitelboim, J. Zanelli,
\emph{Geometry of the $2+1$ black hole,\/} Phys. Rev. D (3), {\bf 48} no 4
(1993), 1506--1525, gr-qc/9302012.

\bibitem{barHLP} 
T. Barbot, \emph{Vari\'et\'es affines radiales de dimension $3$,\/} Bull. Soc. math. France, {\bf 128}
(2000), 347--389.


\bibitem{ba1} 
T. Barbot, \emph{Causal properties of AdS-isometry groups I: 
causal actions and limit sets,\/} math.GT/0509552.



\bibitem{beem} 
J.K. Beem, P.E. Ehrlich, K.L. Easley, 
\emph{Global Lorentzian geometry,\/} Monographs and Textbooks in 
Pure and Applied Mathematics, 2nd ed., {\bf 202}, Marcel Dekker, 
New York, 1996.

\bibitem{BenBon2} 
R. Benedetti, F. Bonsante, 
\emph{Canonical Wick rotations in $3$-dimensional gravity,\/} 
preprint, math.DG/0508485.

\bibitem{brill1} 
D.R. Brill, \emph{Multi-black-hole geometries 
in $(2+1)$-dimensional gravity,\/}  Phys. Rev. D (3), {\bf 53}  (1996),  
no. 8, R4133--R4137, gr-qc/9511022.

\bibitem{brill} 
D.R. Brill, \emph{Black holes and wormholes in $2+1$ 
dimensions,\/}  Mathematical and quantum aspects of relativity 
and cosmology (Pythagoreon, 1998),  143--179, 
Lecture Notes in Phys., {\bf 537}, Springer, Berlin, 2000, gr-qc/9904083.


\bibitem{chrusciel} P. Chrusciel, E. Delay, G.J. Galloway, R. Howard, 
\emph{Regularity of Horizons and the Area Theorems,\/} 
Ann. Henri Poincar\'e, {\bf 2} (2001), 109--178.


\bibitem{gott} S. DeDeo, J.R. Gott III, 
\emph{An eternal time machine in $2+1$ dimensional anti-de Sitter space,\/} 
Phys. Rev. D66 (2002) 084020; Erratum-ibid. D67 (2003) 069902, gr-qc/0212118


\bibitem{bordads} C. Frances, 
\emph{The conformal boundary of anti-de Sitter spacetimes,\/} 
Actes du colloque "AdS/CFT Corespondance", Strasbourg 10-13 Sept. 2003.


\bibitem{fricke} R. Fricke, 
\emph{\"{U}ber die Theorie der automorphen Modulgrupper,\/} 
Nachr. Akad. Wiss. G\"{o}ttingen, (1896), 91--101.

\bibitem{klein} R. Fricke, K. Klein, 
\emph{Vorlesungen des Automorphen Funktionen,\/} 
Teubner, Leipzig, Vol. I (1897), Vol. II (1912).

\bibitem{senovilla} A. Garc\'{\i}a - Parrado, J. M. M. Senovilla, 
\emph{Causal structures and causal boundaries,\/} 
Classical Quantum Gravity, {\bf 22}  (2005),  no. 9, R1--R84, 
gr-qc/0501069.

\bibitem{goldman} W. Goldman, 
\emph{The modular group action on real $SL(2)$-characters of a one-holed 
torus,\/} Geom. Topol., {\bf 7} (2003), 443--486.

\bibitem{taubnut} P. Hajicek, \emph{Extensions of the Taub and NUT 
spaces and extensions of their tangent bundle,\/} 
Commun. Math. Phys., {\bf 17} (1970), 109--126.


\bibitem{matschull} H.J. Matschull, 
\emph{Black hole creation in $2+1$ dimensions,\/} 
Class. Quantum Grav., {\bf 16} (1999), 1069--1095, gr-qc/9809087

\bibitem{mess} G.   Mess, 
\emph{ Lorentz   spacetimes  of   constant   curvature,\/}  Preprint
IHES/M/90/28 (1990).

\bibitem{taubnut2} C.W. Misner, \emph{Taub-NUT space as a counterexample 
to almost anything,\/} Relativity Theory and Astrophysics. AMS Providence,
Rhode Island, Lectures in Appl. Math., {\bf 8} (1967), 160--169.

\bibitem{gravitation} C.W. Misner, K.S. Thorne, J.A. Wheeler, 
\emph{Gravitation,\/} (1973) Freeman, San Fransisco.


\bibitem{Oneillkerr} B. O'Neill, 
\emph{The geometry of Kerr black holes,\/} (1995) A.K. Peters, Wellesley, Massachusetts.

\bibitem{bordpenrose} R. Penrose, 
\emph{Conformal treatment of infinity,\/} 
in Relativit\'e, 
(Lectures, Les Houches, 1963 Summer School of Theoret. Phys., Univ. 
Grenoble, 563--584.


\bibitem{salein} F. Salein, 
\emph{Vari\'et\'es anti-de Sitter de dimension
$3$ poss\'edant un champ de Killing non trivial,\/} C. R. Acad. Sci. Paris
S\'er. I Math., {\bf 324} (1997), 525--530.

\bibitem{townsent} P.K. Townsend, \emph{Black Holes,\/} Lecture notes
of a course given in DAMTP, Cambridge, gr-qc/9707012


\end{thebibliography}
\end{document}